\documentclass[noname,nonblindrev]{informs-noname}

\usepackage{algorithm}
\usepackage{algpseudocode}
\usepackage{amsfonts}
\usepackage{amssymb}
\usepackage{bm}
\usepackage{enumitem}
\usepackage{hyperref}
\usepackage{mathtools}
\usepackage{natbib}
 \bibpunct[, ]{(}{)}{,}{a}{}{,}%
 \def\bibfont{\small}%
\usepackage{pgfplotstable}
\usepackage{relsize}
\usepackage{subcaption}
\usepackage{xcolor}

\DoubleSpacedXI
\TheoremsNumberedThrough
\EquationsNumberedThrough

\hypersetup{colorlinks=true,linkcolor=black,citecolor=black,filecolor=black,urlcolor=black}
\pgfplotsset{compat=1.14}
\usetikzlibrary{pgfplots.groupplots}
\usetikzlibrary{shapes.misc}
\usetikzlibrary{intersections,positioning,fit,calc,shapes.multipart}
\usepgfplotslibrary{fillbetween}

\definecolor{palette1}{RGB}{32,104,170}
\definecolor{palette2}{RGB}{24,96,162}
\definecolor{palette3}{RGB}{77,175,74}
\definecolor{palette4}{RGB}{152,78,163}
\definecolor{palette5}{RGB}{255,127,0}

\begin{document}

\MANUSCRIPTNO{00-0000-0000.00}


\ARTICLEAUTHORS{
    \AUTHOR{Byron Tasseff}
    \AFF{Information Systems and Modeling, Los Alamos National Laboratory, Los Alamos, NM 87545, \EMAIL{btasseff@lanl.gov} \\
         Department of Industrial and Operations Engineering, University of Michigan, Ann Arbor, MI 48109}
    \AUTHOR{Russell Bent}
    \AFF{Applied Mathematics and Plasma Physics, Los Alamos National Laboratory, Los Alamos, NM 87545}
    \AUTHOR{Marina A. Epelman}
    \AFF{Department of Industrial and Operations Engineering, University of Michigan, Ann Arbor, MI 48109}
    \AUTHOR{Donatella Pasqualini}
    \AFF{Information Systems and Modeling, Los Alamos National Laboratory, Los Alamos, NM 87545}
    \AUTHOR{Pascal Van Hentenryck}
    \AFF{H. Milton Stewart School of Industrial and Systems Engineering, Georgia Institute of Technology, Atlanta, GA 30332}}
\RUNAUTHOR{Tasseff et al.}

\RUNTITLE{Exact Mixed-integer Convex Programming Formulation for Optimal Water Network Design}
\TITLE{Exact Mixed-integer Convex Programming Formulation for Optimal Water Network Design}


\ABSTRACT{%
In this paper, we consider the canonical water network design problem, which contains nonconvex potential loss functions and discrete resistance choices with varying costs.
Traditionally, to resolve the nonconvexities of this problem, relaxations of the potential loss constraints have been applied to yield a more tractable mixed-integer convex program (MICP).
However, design solutions to these relaxed problems may not be feasible with respect to the full nonconvex physics.
In this paper, it is shown that, in fact, the original mixed-integer nonconvex program can be reformulated \emph{exactly} as an MICP.
Beginning with a convex program previously used for proving nonlinear network design feasibility, strong duality is invoked to construct a novel, convex primal-dual system embedding all physical constraints.
This convex system is then augmented to form an exact MICP formulation of the original design problem.
Using this novel MICP as a foundation, a global optimization algorithm is developed, leveraging heuristics, outer approximations, and feasibility cutting planes for infeasible designs.
Finally, the algorithm is compared against the previous relaxation-based state of the art in water network design over a number of standard benchmark instances from the literature.}


\maketitle


\section{Introduction}
\label{section:introduction}
This paper considers the problem of optimal water network design, where the layout of the network is known, and the diameter and material of each pipe must be selected from a discrete set to minimize cost while satisfying fixed demand.
The canonical design problem presented throughout the literature excludes common operational components (e.g., pumps), and all demand is assumed to be gravity-fed.
This paper focuses on the development of mathematical programming (nonheuristic) solution techniques for this problem.
However, because of the nonconvexities that appear in the physics of water systems, naive mathematical programming formulations quickly become intractable.
To address this, the problem is typically relaxed via convexification.
Unfortunately, solutions to these relaxed problems may be infeasible with respect to the full nonconvex physics.

To this end, \cite{cherry1951cxvii} introduced a convex program for determining the feasibility of a fixed network governed by linear conservation laws and nonlinear potential loss relationships.
\cite{raghunathan2013global} later revisited this model, using it as a simple tool within a global, relaxation-based design optimization algorithm.
However, the greater potential of a \emph{convex} description for feasibility has gone, for the most part, surprisingly overlooked.
This paper explores this potential, which ultimately provides a new understanding of water network design.
Its contributions include
\begin{itemize}[noitemsep,topsep=0pt]
    \item A convex embedding of \emph{all} physical constraints for a feasible network design;
    \item An \emph{exact} mixed-integer convex reformulation of the optimal design problem;
    \item Development of a new global, relaxation-based algorithm for optimal network design;
    \item A comparison of this algorithm with the previous state of the art \citep{raghunathan2013global};
    \item New lower and upper bounds for open benchmark instances throughout the literature.
\end{itemize}
More broadly, these contributions could be further generalized to optimal network design and expansion planning problems for other important network types exhibiting nonlinear relationships among nodal potentials and flows (e.g., natural gas and crude oil transmission networks).

As background, Section \ref{section:literature_review} reviews optimization techniques for water and other similar network problems, Section \ref{section:problem_formulation} formulates the water network design problem as a mixed-integer nonconvex program, and Section \ref{section:relaxation_based_reformulation} summarizes the contributions of \cite{raghunathan2013global}: (i) their mixed-integer convex programming (MICP) relaxation, (ii) the further relaxation that forms a mixed-integer linear program (MIP), and (iii) the outline of a MIP-based algorithm.
The remaining sections detail this paper's contributions: Section \ref{section:convex_reformulation} derives a convex description of feasible designs, then reformulates the optimal design problem \emph{exactly} as an MICP; Section \ref{section:algorithm} augments Raghunathan's algorithm with novel outer approximations based on the new MICP; Section \ref{section:computational_experiments} compares the new and previous algorithms using instances from the literature; and Section \ref{section:concluding_remarks} concludes the paper.


\section{Literature Review}
\label{section:literature_review}
\emph{Nonlinear networks} refer to a class of networks in which (i) flow is driven by potentials and (ii) potential loss along an edge is a nonlinear function of flow.
Network types with these properties include potable water, natural gas, and crude oil.
Because of their mathematically similar descriptions, optimization methods developed for any of these networks are often easily adapted to the others.
For over fifty years, a variety of techniques have been employed to solve optimization problems that involve these network types \citep{raghunathan2013global}.
\cite{mala2017lost,mala2018lost} provide comprehensive literature reviews of solution techniques used for optimal water network operation and design, both of which are dominated by metaheuristic methods.
\cite{rios2015optimization} provide a similar review of techniques for problems involving natural gas networks.
Finally, \cite{sahebi2014strategic} review methods for optimizing crude oil supply chains.

Outside of mathematical programming, the predominant approaches used for solving water network design problems have been heuristic techniques based on simulation optimization.
Indeed, as \cite{mala2018lost} express, ``... research [has] been trapped, to some extent, in applying new metaheuristic [optimization] methods to relatively simple (from an engineering perspective) design problems, without understanding the principles behind algorithm performance.''
The same review specifies that $84\%$ of the $124$ studies compared use ``stochastic'' methods (e.g., evolutionary and genetic algorithms), $9\%$ use ``deterministic'' methods (e.g., linear and nonlinear programming), and $7\%$ use hybridized methods.
\cite{MAIER2015222} similarly describe the prevalence of metaheuristic optimization techniques in the water resources literature and their associated performance inconsistencies.
For these reasons, as well as due to the lack of optimality guarantees associated with metaheuristic techniques, this paper instead focuses on mathematical programming (specifically, globally optimal mixed-integer nonlinear) approaches for solving water network design problems.

A number of recent studies have developed mathematical programming techniques that have proven to be effective on various nonlinear network problems.
\cite{borraz2016convex} develop a relaxation-based method for natural gas network expansion planning.
\cite{d2015mathematical} provide a survey of methods used throughout water system optimization, which includes both approximation- and relaxation-based techniques for optimal water network operation and design.
\cite{raghunathan2013global} presents a relaxation-based approach for general nonlinear network design and showcases its efficacy on water network design instances.
Our paper considers the algorithm therein to be the state of the art for global nonlinear network design.
As such, their algorithm serves as a foundation for the one developed in this paper, where only the selection of cuts differs.

In conjunction with recent advances in relaxation-based methods for optimizing over nonlinear networks, similar developments have been made in formulating cuts that strengthen these relaxations.
\cite{humpola2013new} develop a number of valid inequalities for nonlinear network design problems and briefly describe their potential for water networks.
\cite{humpola2016valid} extend this work for natural gas.
In both cases, valid cuts are derived from the solution of a nonconvex program and several auxiliary problems.
In contrast, the cuts in this paper are trivially derived from the novel mixed-integer convex reformulation of the original network design problem.

Finally, the notion of a convex feasibility problem for nonlinear networks, originating with \cite{cherry1951cxvii} and exploited computationally by \cite{raghunathan2013global}, is also discussed by \cite{de1994mathematical}.
Here, a convex reformulation of the optimal natural gas transmission problem under restrictive assumptions is presented.
A number of properties of this problem are then described, including solution uniqueness and physical interpretation.
Unlike this paper, however, the problem they consider is not discrete, and their contributions appear to have gone unused throughout the natural gas optimization literature.
Moreover, their work derives the reformulation via variational inequality theory, whereas this paper obtains a more general formulation via Lagrangian duality.


\section{Problem Formulation}
\label{section:problem_formulation}

\subsection{Notation for Sets}
A water distribution network is represented by a directed graph $\mathcal{G} := (\mathcal{N}, \mathcal{A})$, where $\mathcal{N}$ is the set of nodes (i.e., junctions and reservoirs) and $\mathcal{A}$ is the set of arcs (i.e., pipes).
Herein, the set of reservoirs (i.e., source nodes with fixed potentials) is denoted by $\mathcal{S} \subset \mathcal{N}$ and the set of junctions by $\mathcal{J} \subset \mathcal{N}$.
Junctions are modeled as demand nodes (where the demand for flow is nonnegative) and, without loss of generality, $\mathcal{S} \cap \mathcal{J} = \emptyset$.
The set of arcs incident to node $i \in \mathcal{N}$ where $i$ is the tail (respectively, head) of the arc is denoted by $\delta^{+}_{i} := \{(i, j) \in \mathcal{A}\}$ (respectively, $\delta^{-}_{i} := \{(j, i) \in \mathcal{A}\}$).
All arcs incident to a source node $i \in \mathcal{S}$ are assumed to be outgoing, i.e., $\delta_{i}^{-} = \emptyset$ and $\delta_{i}^{+} \neq \emptyset$, and arcs with tails at demand nodes $i \in \mathcal{J}$ have heads also at demand nodes.
Finally, the design problem of Section \ref{subsection:minlp} involves selecting from a set of resistances $\mathcal{R}_{a} := \left\{p_{1}, p_{2}, \dots, p_{\lvert\mathcal{R}_{a}\rvert}\right\}$ for each $a \in \mathcal{A}$.

\subsection{Physical Feasibility}
This section describes the variables and constraints required to model the physics of gravity-fed water networks given a \emph{fixed} selection of pipe resistances $r_{a}$, $a \in \mathcal{A}$.
In the constraints that follow, $q_{a}$, $a \in \mathcal{A}$, denote variables representing the flow of water across each arc (expressed as a volumetric flow rate in $\textnormal{m}^{3}/\textnormal{s}$).
Nodal potentials are denoted by the variables $h_{i}$, $i \in \mathcal{N}$, where each represents the total hydraulic head in units of length ($\textnormal{m}$).
The total hydraulic head (hereafter referred to as ``head'') assimilates elevation and pressure heads at a node, while the velocity head is neglected.

\paragraph{Flow Bounds}
When $q_{a}$ is positive (negative), flow on arc $a := (i, j)$ travels from node $i$ to $j$ ($j$ to $i$).
Flow is often bounded by physical capacity, engineering judgment, or network analysis.
Herein,
\begin{equation}
	\underline{q}_{a} \leq q_{a} \leq \overline{q}_{a}, \; \forall a \in \mathcal{A} \label{eqn:feas-flow-rate-bounds}.
\end{equation}
For example, the maximum speed of flow along arc $a \in \mathcal{A}$, $\overline{v}_{a}$, is often used to estimate the bounds $\underline{q}_{a} = -\frac{\pi}{4} \overline{v}_{a} D_{a}^{2}$ and $\overline{q}_{a} = \frac{\pi}{4} \overline{v}_{a} D_{a}^{2}$, where $D_{a}$ is the fixed diameter of the pipe indexed by $a \in \mathcal{A}$.

\paragraph{Head Bounds}
For each reservoir $i \in \mathcal{S}$, the head is fixed at a constant value $h_{i}^{s}$, i.e.,
\begin{equation}
	h_{i} = h_{i}^{s}, \; \forall i \in \mathcal{S} \label{eqn:feas-fixed-head-at-source}.
\end{equation}
For each junction $i \in \mathcal{J}$, a predefined minimum head $\underline{h}_{i}$ must be satisfied.
Upper bounds on heads can also be provided or implied by network data.
For example, this paper assumes
\begin{equation}
	\underline{h}_{i} \leq h_{i} \leq \overline{h}_{i} = \max_{j \in \mathcal{S}}\{h_{j}^{s}\}, \; \forall i \in \mathcal{J} \label{eqn:feas-head-bounds}.
\end{equation}

\paragraph{Conservation of Flow at Demand Nodes}
Flow must be delivered throughout the network in order to satisfy fixed nonnegative demand, $d_{i}$, at all demand nodes $i \in \mathcal{J}$. That is,
\begin{equation}
	\sum_{\mathclap{a \in \delta^{-}_{i}}} q_{a} - \sum_{\mathclap{a \in \delta^{+}_{i}}} q_{a} = d_{i}, \; \forall i \in \mathcal{J} \label{eqn:feas-flow-conservation}.
\end{equation}

\paragraph{Head Loss Relationships}
In water networks, flow along an arc is induced by the difference in head between the two nodes connected by that arc.
The relationships that link flow and head are commonly referred to as the ``head loss equations'' and are generally of the form
\begin{equation}
	h_{i} - h_{j} = \phi_{a}(q_{a}), ~ \forall a := (i, j) \in \mathcal{A},
\end{equation}
where $\phi_{a} : \mathbb{R} \to \mathbb{R}$ is a strictly increasing function with rotational symmetry about the origin.
The most common head loss relationships include the Darcy-Weisbach equation,
\begin{equation}
	h_{i} - h_{j} = \frac{8 L_{a} \tau_{a} q_{a} \lvert q_{a} \rvert}{\pi^{2} g D_{a}^{5}}, \label{eqn:darcy-weisbach}
\end{equation}
and the Hazen-Williams equation (where the constant $10.7$ is in standard units),
\begin{equation}
    h_{i} - h_{j} = \frac{10.7 L_{a} q_{a} \lvert q_{a} \rvert^{0.852}}{\kappa_{a}^{1.852} D_{a}^{4.8704}} \label{eqn:hazen-williams}.
\end{equation}
Here, $L_{a}$ is the pipe length, $\tau_{a}$ is the friction factor, $g$ is gravitational acceleration, and $\kappa_{a}$ is the roughness, which depends on the pipe material.
In Equation \eqref{eqn:darcy-weisbach}, $\tau_{a}$ depends on $q_{a}$ in a nonlinear manner.
However, in the mathematical programming literature, $\tau_{a}$ is often fixed to a constant, which removes the term's nonlinearity in $q_{a}$ \citep{gleixner2012towards,verleye2013optimising}.

When all terms \emph{except} $h_{i}$, $h_{j}$, and $q_{a}$ are fixed, both head loss equations reduce to
\begin{equation}
	h_{i} - h_{j} = L_{a} r_{a} q_{a} \lvert q_{a} \rvert^{\alpha - 1}, \; \forall a := (i, j) \in \mathcal{A} \label{eqn:feas-head-loss}.
\end{equation}
Here, $\alpha$ denotes the exponent required by Equation \eqref{eqn:darcy-weisbach} or \eqref{eqn:hazen-williams} (i.e., $2$ or $1.852$, respectively), and $r_{a}$, $a \in \mathcal{A}$, denotes the resistance per unit length.
The resistance per unit length comprises all non-length constant terms appearing in Equations \eqref{eqn:darcy-weisbach} and \eqref{eqn:hazen-williams} and is in units of $(\textrm{m}^3 / \textrm{s})^{-\alpha}$.

For fixed resistances $r$, the nonconvex formulation for water network feasibility is thus
\begin{equation} \tag{NLP$(r)$}
\begin{aligned}
	& \textnormal{Physical bounds: Constraints} ~ \eqref{eqn:feas-flow-rate-bounds}, \eqref{eqn:feas-fixed-head-at-source}, \eqref{eqn:feas-head-bounds} \\
	& \textnormal{Flow conservation: Constraints} ~ \eqref{eqn:feas-flow-conservation} \\
	& \textnormal{Head loss relationships: Constraints} ~ \eqref{eqn:feas-head-loss}.
\end{aligned} \label{eqn:nlp}
\end{equation}

\subsection{Optimal Network Design}
\label{subsection:minlp}
To formulate the design problem, \eqref{eqn:nlp} is coupled with the combinatorial problem of selecting one resistance from a predefined set of resistances for each arc,  $\mathcal{R}_{a}$ for $a \in \mathcal{A}$, while minimizing the overall cost of network design.
To model the disjunction representing discrete resistance choices, each $q_{a}$ is first decomposed into $\lvert \mathcal{R}_{a} \rvert$ binary variables $x_{ap}$ and continuous variables $q_{ap}$, i.e.,
\begin{equation}
    q_{a} := \sum_{\mathclap{p \in \mathcal{R}_{a}}} q_{ap}, ~ \forall a \in \mathcal{A} \label{eqn:decomposition-of-q},
\end{equation}
\begin{equation}
    \underline{q}_{ap} x_{ap} \leq q_{ap} \leq \overline{q}_{ap} x_{ap}, ~ x_{ap} \in \mathbb{B}, ~ \forall a := (i, j) \in \mathcal{A}, ~ \forall p \in \mathcal{R}_{a} \label{eqn:minlp-flow-rate-bounds}.
\end{equation}
Here, $x_{ap} = 1$ when $p \in \mathcal{R}_{a}$ is selected and zero otherwise.
From Constraint \eqref{eqn:minlp-flow-rate-bounds}, it follows that $q_{ap}$ is nonzero only when $x_{ap} = 1$.
Also note that $\underline{q}_{ap} = -\frac{\pi}{4} \overline{v}_{a} D_{ap}^{2}$ and $\overline{q}_{ap} = \frac{\pi}{4} \overline{v}_{a} D_{ap}^{2}$ are often used as bounds, as each $p \in \mathcal{R}_{a}$ is typically derived from a unique pipe diameter $D_{ap}$.
Since only one resistance may be selected per pipe indexed by $a \in \mathcal{A}$, we include the additional constraints
\begin{equation}
    \sum_{\mathclap{p \in \mathcal{R}_{a}}} x_{ap} = 1, ~ \forall a \in \mathcal{A} \label{eqn:ne-resistance-selection}.
\end{equation}
The head loss Constraints \eqref{eqn:feas-head-loss} are then expanded to formulate the constraints
\begin{equation}
    h_{i} - h_{j} = L_{a} \sum_{\mathclap{p \in \mathcal{R}_{a}}} p q_{ap} \lvert q_{ap} \rvert^{\alpha - 1}, ~ \forall a := (i, j) \in \mathcal{A} \label{eqn:minlp-head-loss}.
\end{equation}
Finally, the objective function, $\eta(x)$, for the optimal design problem is written as
\begin{equation}
	\eta(x) = \sum_{\mathclap{a \in \mathcal{A}}} L_{a} \sum_{\mathclap{p \in \mathcal{R}_{a}}} c_{ap} x_{ap} \label{eqn:ne-obj},
\end{equation}
where $c_{ap}$ is the cost per unit length of installing a pipe along $a \in \mathcal{A}$ with resistance $p \in \mathcal{R}_{a}$.

These modifications allow the optimal water network design problem to be written as
\begin{equation}\tag{MINLP}\begin{aligned}
    & \text{minimize}
	& & \textnormal{Objective function:} ~ \eta(x) ~ \textnormal{of Equation} ~ \eqref{eqn:ne-obj} \\
    & \text{subject to}
	& & \textnormal{Physical bounds: Constraints} ~ \eqref{eqn:feas-fixed-head-at-source}, \eqref{eqn:feas-head-bounds}, \eqref{eqn:minlp-flow-rate-bounds} \\
	& & & \textnormal{Flow conservation: Constraints} ~ \eqref{eqn:feas-flow-conservation} \\
	& & & \textnormal{Resistance selection: Constraints} ~ \eqref{eqn:ne-resistance-selection} \\
	& & & \textnormal{Head loss relationships: Constraints} ~ \eqref{eqn:minlp-head-loss}.
\end{aligned}\label{eqn:minlp}\end{equation}
Here, Constraints \eqref{eqn:feas-flow-conservation} employ the definitions of $q_{a}$ described in Equations \eqref{eqn:decomposition-of-q}.
Note that \eqref{eqn:minlp} is mixed-integer \emph{nonconvex} because of the nonconvex Constraints \eqref{eqn:minlp-head-loss}.
\cite{raghunathan2013global} addresses this challenge via a convex relaxation of the complicating constraints, followed by linearization of the resulting MICP.
Their method is reviewed in Section \ref{section:relaxation_based_reformulation}.
This paper, on the other hand, addresses this challenge via an \emph{exact} convex reformulation of \eqref{eqn:nlp} to describe feasibility, followed by an \emph{exact} MICP reformulation of \eqref{eqn:minlp}.
These are further described in Section \ref{section:convex_reformulation}.


\section{Relaxation-based Reformulation and Algorithm}
\label{section:relaxation_based_reformulation}
This section reviews the methods of \cite{raghunathan2013global}, which this paper uses as its foundations.
First, their relaxed MICP of the design problem is presented in Section \ref{subsec:micp-relaxation}.
Then, an exact MIP reformulation of \eqref{eqn:minlp}, based on an outer approximation of the MICP, is presented in Section \ref{subsec:mip-relaxation}.
A convex program for determining the feasibility of a design is presented in Section \ref{subsec:cnlp}.
Finally, a simple global algorithm, which leverages the relaxed MIP formulation described in Section \ref{subsec:mip-relaxation}, is presented in Section \ref{subsec:milpr-algorithm}.
These details prepare the reader for \emph{this} study's \emph{exact} MICP reformulation described in Section \ref{section:convex_reformulation}, as well as for the relaxation-based global algorithm presented in Section \ref{section:algorithm}.

\subsection{Mixed-integer Convex Relaxation of (MINLP)}
\label{subsec:micp-relaxation}
This section relaxes the nonconvex head loss constraints of \eqref{eqn:minlp} via an outer convexification to form a relaxed MICP of the original optimal design problem.
This is accomplished by partitioning Constraints \eqref{eqn:minlp-head-loss} into their symmetric positive and negative components.
To begin, variables $q_{ap}^{\pm}$ are introduced, denoting nonnegative flows in the two directions along arc $a \in \mathcal{A}$, with
\begin{equation}
    q_{a} := \sum_{\mathclap{p \in \mathcal{R}_{a}}} \left(q^{+}_{ap} - q^{-}_{ap}\right), ~ \forall a \in \mathcal{A} \label{eqn:flow-equations}
\end{equation}
replacing Equations \eqref{eqn:decomposition-of-q}.
Next, the bound Constraints \eqref{eqn:minlp-flow-rate-bounds} in \eqref{eqn:minlp} are rewritten as
\begin{equation}
    0 \leq q^{\pm}_{ap} \leq \overline{q}^{\pm}_{ap} x_{ap}, ~ x_{ap} \in \mathbb{B}, ~ \forall a \in \mathcal{A}, ~ \forall p \in \mathcal{R}_{a} \label{eqn:micpr-flow-bounds},
\end{equation}
where $\overline{q}_{ap}^{+} = \max\{0, \overline{q}_{ap}\}$ and $\overline{q}_{ap}^{-} = \max\{0, -\underline{q}_{ap}\}$ replace the bounds of Constraints \eqref{eqn:minlp-flow-rate-bounds}.
Nonnegative head difference variables $\Delta h_{a}^{\pm}$ are similarly introduced to denote head loss in the two possible flow directions.
These are related to the original head variables $h_{i}$, $i \in \mathcal{N}$, via the constraints
\begin{equation}
    \Delta h_{a}^{+} - \Delta h_{a}^{-} = h_{i} - h_{j}, ~ \forall a := (i,j) \in \mathcal{A} \label{eqn:micpr-head-equality}.
\end{equation}

Next, binary variables $y_{a} \in \mathbb{B}$, $a \in \mathcal{A}$, are introduced to denote the direction of flow along each arc, where, for $a=(i,j)$, $y_{a} = 1$ implies flow from $i$ to $j$ and $y_{a} = 0$ implies flow from $j$ to $i$, i.e.,
\begin{subequations}\label{eqn:micpr-direction-bounds}\begin{gather}
	0 \leq q_{ap}^{+} \leq \overline{q}_{ap}^{+} y_{a}, ~ 0 \leq q_{ap}^{-} \leq \overline{q}_{ap}^{-} (1 - y_{a}), ~ y_{a} \in \mathbb{B}, ~ \forall a \in \mathcal{A}, ~ \forall p \in \mathcal{R}_{a} \\
	0 \leq \Delta h_{a}^{+} \leq \Delta \overline{h}_{a}^{+} y_{a}, ~ 0 \leq \Delta h_{a}^{-} \leq \Delta \overline{h}_{a}^{-} (1 - y_{a}), ~ y_{a} \in \mathbb{B}, ~ \forall a \in \mathcal{A}.
\end{gather}\end{subequations}
Here, each head difference bound $\Delta \overline{h}_{a}^{\pm}$ is derived from the lower and upper bounds on $h$ in Constraints \eqref{eqn:feas-head-bounds}.
Next, the right-hand sides in Equations \eqref{eqn:feas-head-loss} are decomposed into two convex functions representing head loss in the positive and negative directions.
This decomposition implies
\begin{equation}
    \Delta h_{a}^{\pm} = L_{a} \sum_{\mathclap{p \in \mathcal{R}_{a}}} p \left(q_{ap}^{\pm}\right)^{\alpha}, ~ \forall a \in \mathcal{A} \label{eqn:head-loss-decomposed}.
\end{equation}
Recalling that $q_{ap}^{\pm}$'s can be nonzero for only one $p \in \mathcal{R}_{a}$, Equations \eqref{eqn:head-loss-decomposed} are then relaxed as
\begin{equation}
    L_{a} p \left(q_{ap}^{\pm}\right)^{\alpha} \leq \Delta h_{a}^{\pm}, ~ \forall a \in \mathcal{A}, ~ \forall p \in \mathcal{R}_{a} \label{eqn:micpr-head-loss}.
\end{equation}
Note that $(0, 0)$ and $\left(\overline{q}_{ap}^{\pm}, L_{a} p \left(\overline{q}_{ap}^{\pm}\right)^{\alpha}\right)$ are endpoints of the lines that upper-bound the strictly convex right-hand terms of Equations \eqref{eqn:head-loss-decomposed}, with the slopes of these lines calculated as  $\frac{L_{a} p \left(\overline{q}_{ap}^{\pm}\right)^{\alpha} - 0}{\overline{q}_{ap}^{\pm} - 0} = L_{a} p \left(\overline{q}_{ap}^{\pm}\right)^{\alpha - 1}$.
The convex relaxations in Constraints \eqref{eqn:micpr-head-loss} are then linearly upper-bounded via
\begin{equation}
    \Delta h_{a}^{\pm} \leq L_{a} \sum_{\mathclap{p \in \mathcal{R}_{a}}} \left[p \left(\overline{q}_{ap}^{\pm}\right)^{\alpha - 1} q_{ap}^{\pm}\right], ~ \forall a \in \mathcal{A} \label{eqn:micpr-head-difference-ub}.
\end{equation}

These constraints give rise to an MICP \emph{relaxation} of \eqref{eqn:minlp}, namely,
\begin{equation}\tag{MICP-R}\begin{aligned}
    & \text{minimize}
	& & \textnormal{Objective function:} ~ \eta(x) ~ \textnormal{of Equation} ~ \eqref{eqn:ne-obj} \\
    & \text{subject to}
	& & \textnormal{Physical bounds: Constraints} ~ \eqref{eqn:feas-fixed-head-at-source}, \eqref{eqn:feas-head-bounds}, \eqref{eqn:micpr-flow-bounds} \\
	& & & \textnormal{Flow conservation: Constraints} ~ \eqref{eqn:feas-flow-conservation} \\
	& & & \textnormal{Resistance selection: Constraints} ~ \eqref{eqn:ne-resistance-selection} \\
	& & & \textnormal{Head difference relationships: Constraints} ~ \eqref{eqn:micpr-head-equality}, \eqref{eqn:micpr-head-loss}, \eqref{eqn:micpr-head-difference-ub} \\
    & & & \textnormal{Direction-related inequalities: Constraints} ~ \eqref{eqn:micpr-direction-bounds}.
\end{aligned}\label{eqn:micpr}\end{equation}
Here, Constraints \eqref{eqn:feas-flow-conservation} employ the definitions of $q_{a}$ described in Equations \eqref{eqn:flow-equations}.
The validity of this relaxed formulation was proven by \cite{raghunathan2013global}.
However, note that because \eqref{eqn:micpr} is a relaxation of \eqref{eqn:minlp}, network design solutions that are feasible to \eqref{eqn:micpr} may not be feasible with respect to the original, nonconvex head loss constraints of \eqref{eqn:minlp}.

\subsection{Mixed-integer Linear Reformulation of (MINLP)}
\label{subsec:mip-relaxation}
\cite{raghunathan2013global} solves \eqref{eqn:minlp} via a global, relaxation-based MIP algorithm.
The algorithm leverages much of \eqref{eqn:micpr}, linear outer approximations of Constraints \eqref{eqn:micpr-head-loss}, and linear feasibility cuts for integer solutions $\hat{x}$ that are infeasible to \eqref{eqn:minlp}.
This section restates these cuts, then develops a MIP reformulation of \eqref{eqn:minlp}, representing Raghunathan's theoretical contributions.

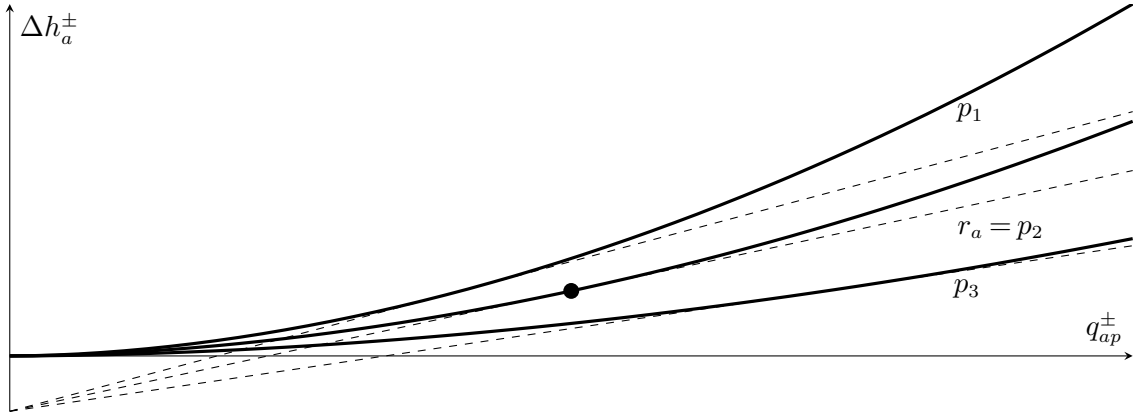
\begin{figure}[t]
    \begin{tikzpicture}[declare function={f_1(\x)=ifthenelse(\x>=0,100*0.422657*\x*abs(\x)^(0.852),100*0.422657*\x*(-\x)^(0.852));},declare function={f_2(\x)=ifthenelse(\x>=0,50*0.422657*\x*abs(\x)^(0.852),50*0.422657*\x*(-\x)^(0.852));}, declare function={f_3(\x)=ifthenelse(\x>=0,150*0.422657*\x*abs(\x)^(0.852),150*0.422657*\x*(-\x)^(0.852));}]
\begin{axis}[
	axis on top,
	legend pos = outer north east,
	axis lines = center,
	xtick style={draw=none},
	ytick style={draw=none},
	xticklabels={,,},
	yticklabels={,,},
	height=7.0cm,
	width=\columnwidth,
	xmin=0.0,
	xmax=0.2,
	xlabel = $q_{ap}^{\pm}$,
       ylabel = $\Delta h_{a}^{\pm}$,
	legend style = {cells={align=left}},
	legend cell align = {left}]
	\addplot[black,very thick,samples=50,domain=0.0:0.2,name path=f] {f_1(x)} node at (0.1765,1.15) {$r_{a} = p_{2}$};
	\addplot[black,very thick,samples=50,domain=0.0:0.2,name path=f2] {f_2(x)} node[below,pos=0.75] {$p_{3}$};
	\addplot[black,very thick,samples=50,domain=0.0:0.2,name path=f3] {f_3(x)} node[below,pos=0.75] {$p_{1}$};
	\addplot[mark=none,dashed,black] coordinates {(0.0,-0.5063) (0.2,1.6949)};
	\addplot[mark=none,dashed,black] coordinates {(0.0,-0.5063) (0.2,1.0076)};
	\addplot[mark=none,dashed,black] coordinates {(0.0,-0.5063) (0.2,2.2336)};
	\fill (0.1,0.5943) circle[radius=3pt];
\end{axis}
\end{tikzpicture}
    \caption{Depiction of equal intercept outer approximations (dashed) of convexified head loss relations (solid) for a hypothetical case where $\lvert \mathcal{R}_{a} \rvert = 3$, $r_{a} = p_{2}$, and where $\tilde{q}_{ar_{a}}$ is depicted by the single point.}
    \label{fig:equal-intercept}
\end{figure}

\paragraph{Outer Approximation Cutting Planes}
Since Constraints \eqref{eqn:micpr-head-loss} are convex, they are easily linearized via outer approximation.
However, instead of applying traditional first-order outer approximations for \emph{each} Constraint \eqref{eqn:micpr-head-loss}, \citet{raghunathan2013global} derives aggregate outer approximations based on the notion of outer-approximating lines with \emph{equal intercepts} for all $p \in \mathcal{R}_{a}$.
This is illustrated in Figure \ref{fig:equal-intercept} for an instance where $\lvert \mathcal{R}_{a} \rvert = 3$ and $r_{a} = p_{2}$.
Note that a standard outer approximation of Constraint \eqref{eqn:micpr-head-loss} for $p \in \mathcal{R}_{a}$, based on the first-order Taylor expansion of its left-hand side at $\tilde{q}_{ap}^{\pm}$, is
\begin{equation}
    \left(1 - \alpha\right) L_{a} p \left(\tilde{q}_{ap}^{\pm}\right)^{\alpha} + \alpha L_{a} p \left(\tilde{q}_{ap}^{\pm}\right)^{\alpha - 1} q^{\pm}_{ap} \leq \Delta h_{a}^{\pm}. \label{eqn:imip-oa-generic}
\end{equation}
Thus, if we fix a point $\tilde{q}_{ar_{a}}^{\pm} \in [0, \overline{q}_{ar_{a}}^{\pm}]$ for some chosen $r_{a} \in \mathcal{R}_{a}$ and $a \in \mathcal{A}$, for the remaining $p \in \mathcal{R}_{a}$, the outer approximations will have the \emph{same constant intercept} for $\tilde{q}_{ap}^{\pm} := \tilde{q}_{ar_{a}}^{\pm}(r_{a}  / p)^{1 / \alpha}$.

Recall that variables $q_{ap}^{\pm}$ become nonzero only when the corresponding resistance $p$ is active (i.e., $x_{ap} = 1$).
Also note that, since only one flow direction is chosen per arc in \eqref{eqn:micpr}, cuts can be strengthened through multiplication of constants with affine expressions of $y_{a}$.
Exploiting these properties, aggregating over all $r_{a} \in \mathcal{R}_{a}$ and $\tilde{q}_{ar_{a}}^{\pm}$ allows the rewriting of Constraints \eqref{eqn:micpr-head-loss} as
\begin{subequations}\begin{align}
    \tilde{\tau}_{ar_{a}}^{+} y_{a} + \alpha \sum_{\mathclap{p \in \mathcal{R}_{a}}} p \left(\tilde{q}_{ap}^{+}\right)^{\alpha - 1} q_{ap}^{+}
    &\leq \frac{\Delta h_{a}^{+}}{L_{a}}, \, \forall a \in \mathcal{A}, \, \forall r_{a} \in \mathcal{R}_{a}, \, \forall \tilde{q}_{ar_{a}}^{+} \in \mathcal{Q}^{+}_{ar_{a}} \\
    \tilde{\tau}_{ar_{a}}^{-} \left(1 - y_{a}\right) + \alpha \sum_{\mathclap{p \in \mathcal{R}_{a}}} p \left(\tilde{q}_{ap}^{-}\right)^{\alpha - 1} q_{ap}^{-}
    &\leq \frac{\Delta h_{a}^{-}}{L_{a}}, \, \forall a \in \mathcal{A}, \, \forall r_{a} \in \mathcal{R}_{a}, \, \forall \tilde{q}_{ar_{a}}^{-} \in \mathcal{Q}^{-}_{ar_{a}},
\end{align}\label{eqn:imip-head-loss-aggregated-4}\end{subequations}
where $\tilde{\tau}_{ar_{a}}^{\pm} := \left(1 - \alpha\right) r_{a} \left(\tilde{q}_{ar_{a}}^{\pm}\right)^{\alpha}$ for brevity, and $\mathcal{Q}^{\pm}_{ar_{a}} = [0, \overline{q}_{ar_{a}}^{\pm}]$.
\citet{raghunathan2013global} shows that each cut of Constraints \eqref{eqn:imip-head-loss-aggregated-4} is, in fact, \emph{stronger} than the set of standard disaggregated cuts that outer-approximate Constraints \eqref{eqn:micpr-head-loss}.
Also, as shown in Section \ref{subsec:milpr-algorithm}, algorithmically, $\mathcal{Q}^{\pm}_{ar_{a}}$ can instead be replaced by the finite sets $\tilde{\mathcal{Q}}^{\pm}_{ar_{a}}$ to linearly \emph{approximate} \eqref{eqn:micpr} rather than reproduce it.

\paragraph{Feasibility Cutting Planes}
Since \eqref{eqn:micpr} is a relaxation of \eqref{eqn:minlp}, a design solution to \eqref{eqn:micpr} is not guaranteed to be physically feasible.
To address this, let $\bar{\mathcal{X}}$ denote the set of designs represented by binary vectors $\bar{x}$ satisfying Constraints \eqref{eqn:ne-resistance-selection} that are \emph{not} physically feasible.
Then, any infeasible designs permitted by \eqref{eqn:micpr} will be excluded by the set of linear cuts
\begin{equation}
    \mathlarger{\mathlarger{\sum}}_{a \in \mathcal{A}} \left[\left(\sum_{p \in \mathcal{R}_{a} : \bar{x}_{ap} = 1} x_{ap}\right) - \left(\sum_{p \in \mathcal{R}_{a} : \bar{x}_{ap} = 0} x_{ap}\right)\right] \leq \lvert \mathcal{A} \rvert - 1, ~ \forall \bar{x} \in \bar{\mathcal{X}} \label{eqn:no-good-cut}.
\end{equation}
Each Constraint \eqref{eqn:no-good-cut} is a traditional combinatorial ``no good'' cut, which removes \emph{one} combination of resistances (i.e., one network design) from the space of solutions feasible to \eqref{eqn:micpr}.

\paragraph{Mixed-integer Linear Reformulation}
Combining much of \eqref{eqn:micpr} with the cuts described in this section, the infinite MIP reformulation of \eqref{eqn:minlp}, which abstractly describes the theoretical contributions of and foundational formulation used by \cite{raghunathan2013global}, is
\begin{equation}\tag{MIP-R}\begin{aligned}
    & \text{minimize}
	& & \textnormal{Objective function:} ~ \eta(x) ~ \textnormal{of Equation} ~ \eqref{eqn:ne-obj} \\
    & \text{subject to}
	& & \textnormal{Physical bounds: Constraints} ~ \eqref{eqn:feas-fixed-head-at-source}, \eqref{eqn:feas-head-bounds}, \eqref{eqn:micpr-flow-bounds} \\
	& & & \textnormal{Flow conservation: Constraints} ~ \eqref{eqn:feas-flow-conservation} \\
	& & & \textnormal{Resistance selection: Constraints} ~ \eqref{eqn:ne-resistance-selection} \\
	& & & \textnormal{Head difference relationships: Constraints} ~ \eqref{eqn:micpr-head-equality}, \eqref{eqn:micpr-head-difference-ub}, \eqref{eqn:imip-head-loss-aggregated-4} \\
    & & & \textnormal{Direction-related inequalities: Constraints} ~ \eqref{eqn:micpr-direction-bounds} \\
    & & & \textnormal{Feasibility cutting planes: Constraints} ~ \eqref{eqn:no-good-cut}.
\end{aligned}\label{eqn:milpr}\end{equation}
Algorithmically, the sets $\mathcal{Q}^{\pm}_{ar_{a}}$ and $\bar{\mathcal{X}}$ of Constraints \eqref{eqn:imip-head-loss-aggregated-4} and \eqref{eqn:no-good-cut} are replaced by the initially empty finite sets, $\tilde{\mathcal{Q}}^{\pm}_{ar_{a}}$ and $\tilde{\mathcal{X}}$, and progressively augmented during the algorithm's execution.
Hereafter, we refer to the relaxation of \eqref{eqn:milpr} as (MIP-RR).
To refine this relaxation, the convex oracle algorithmically used for determining whether a design is contained in $\bar{\mathcal{X}}$ is described in Section \ref{subsec:cnlp}.

\subsection{Convex Method for Determining Design Feasibility}
\label{subsec:cnlp}
This section reviews the work of \cite{cherry1951cxvii} and \cite{raghunathan2013global}, who provide a convex program for determining the feasibility of a nonlinear network with fixed resistances $r$.
These resistances correspond to an integer solution $\hat{x}$ satisfying Constraints \eqref{eqn:ne-resistance-selection} via the relationship
\begin{equation}
    r_{a} \in \{p \in \mathcal{R}_{a} : \hat{x}_{ap} = 1\}, ~ \forall a \in \mathcal{A} \label{eqn:x-relation}.
\end{equation}
Thus, this feasibility-testing method can be used to determine whether $\hat{x}$ is contained in $\bar{\mathcal{X}}$.
Letting $q_{a}^{\pm}$ denote nonnegative directed flow variables along $a \in \mathcal{A}$ and using the definition $q_{a} := q_{a}^{+} - q_{a}^{-}$, the studies ultimately propose a strictly convex programming problem similar to
\begin{equation}\tag{P$(r)$}\begin{aligned}
    & \underset{q^{\pm} \geq 0}{\text{minimize}}
    & & \sum_{a \in \mathcal{A}} \frac{L_{a} r_{a}}{1 + \alpha} \left[(q_{a}^{+})^{1 + \alpha} + (q_{a}^{-})^{1 + \alpha}\right] -
    \sum_{i \in \mathcal{S}} h_{i}^{s} \sum_{\mathclap{a \in \delta^{+}_{i}}} q_{a} \\
    & \text{subject to}
    & & \sum_{\mathclap{a \in \delta^{-}_{i}}} q_{a} - \sum_{\mathclap{a \in \delta^{+}_{i}}} q_{a} = d_{i}, \; \forall i \in \mathcal{J}.
\end{aligned}\label{eqn:p}\end{equation}

\cite{raghunathan2013global} proves three important properties of \eqref{eqn:p}: (i) its solution $\hat{q}$ is unique; (ii) the dual multipliers of flow conservation constraints correspond to unique heads $\hat{h}$; and (iii) the primal-dual solution $(\hat{q}, \hat{h})$ satisfies flow conservation and head loss Constraints \eqref{eqn:feas-flow-conservation} and \eqref{eqn:feas-head-loss}.
The ability to find a solution $(\hat{q}, \hat{h})$ satisfying these physical equations by solving \eqref{eqn:p} is appealing because the feasibility of $r$ can be tested by checking whether $(\hat{q}, \hat{h})$ satisfies Constraints \eqref{eqn:feas-flow-rate-bounds} and \eqref{eqn:feas-head-bounds}.
Raghunathan uses this method to guarantee global convergence of their relaxation-based algorithm.

\subsection{Global Optimization Algorithm}
\label{subsec:milpr-algorithm}
This section outlines a simple algorithm for solving \eqref{eqn:minlp} to global optimality via the iterative solution and augmentation of (MIP-RR).
Algorithm \ref{alg:global-milpr} leverages the outer approximation and feasibility cuts established in Section \ref{subsec:mip-relaxation}.
In Line \ref{line:g-init-sets}, the outer approximation point sets and infeasible design set, $\tilde{\mathcal{Q}}^{\pm}$ and $\tilde{\mathcal{X}}$, are initialized as empty.
In Line \ref{line:g-solve-milpr-1}, the relaxed problem is solved, and the solution components $(\hat{q}, \hat{h}, \hat{x})$ are stored.
In Line \ref{line:g-check-feasibility}, feasibility of the design $\hat{x}$ is determined (e.g., via the method of Section \ref{subsec:cnlp}).
If the design is found to be physically infeasible, outer approximations of head loss constraints are added in Line \ref{line:g-add-oa}.
A feasibility cut that excludes $\hat{x}$ is then added in Line \ref{line:g-add-feas}.
Finally, a solution to the new relaxed problem, with the aforementioned cuts, is obtained in Line \ref{line:g-solve-milpr-2}.
These steps are repeated until a feasible solution $\hat{x}$ to \eqref{eqn:minlp} is identified.
Since this $\hat{x}$ is discovered via the solution of sequential \emph{relaxations} of \eqref{eqn:minlp} and the cuts do not exclude feasible solutions to \eqref{eqn:minlp}, the design obtained is guaranteed to be globally optimal.

\begin{algorithm}[!ht]
    \caption{A MIP relaxation-based global optimization algorithm for \eqref{eqn:minlp}.}
    \label{alg:global-milpr}
    \begin{algorithmic}[1]
        \State{$\tilde{\mathcal{Q}}_{ap}^{\pm} \gets \emptyset, ~ \forall a \in \mathcal{A}, ~ \forall p \in \mathcal{R}_{a}$; $\tilde{\mathcal{X}} \gets \emptyset$. \label{line:g-init-sets}}
        \State{$(\hat{q}, \hat{h}, \hat{x})$ $\gets$ Solve (MIP-RR). \label{line:g-solve-milpr-1}}
        \While{$\hat{x}$ is infeasible to \eqref{eqn:minlp} \label{line:g-check-feasibility}}
            \State{$\tilde{\mathcal{Q}}_{ap}^{\pm} \gets \tilde{\mathcal{Q}}_{ap}^{\pm} \cup \{\hat{q}_{ap}^{\pm}\}, ~ \forall a \in \mathcal{A}, ~ \forall p \in \mathcal{R}_{a}$. \label{line:g-add-oa}}
            \State{$\tilde{\mathcal{X}} \gets \tilde{\mathcal{X}} \cup \{\hat{x}\}$. \label{line:g-add-feas}}
            \State{$(\hat{q}, \hat{h}, \hat{x})$ $\gets$ Solve (MIP-RR). \label{line:g-solve-milpr-2}}
        \EndWhile
    \end{algorithmic}
\end{algorithm}

In practice, the algorithm developed by \cite{raghunathan2013global} is more sophisticated than Algorithm \ref{alg:global-milpr}.
Specifically, it exploits the linearization-based linear programming/nonlinear programming branch and bound (LP/NLP-BB) framework developed by \cite{quesada1992lp}.
This permits more flexibility than Algorithm \ref{alg:global-milpr} through the use of ``callback'' features available in many MIP solvers.
First, outer approximations are added in other parts of the search tree, not just integer solutions to (MIP-RR), as in Line \ref{line:g-add-oa}.
Moreover, outer approximation points $\tilde{\mathcal{Q}}^{\pm}$ are more thoughtfully selected.
Second, Raghunathan develops heuristic procedures, internal to the search, that can recover feasible solutions from fractional solutions and integer solutions \emph{infeasible} to \eqref{eqn:minlp}.

In this paper, the majority of Raghunathan's contributions are used to devise a new algorithm, which is further discussed in Section \ref{section:algorithm}.
The two algorithms primarily differ in the selection of outer approximation cuts.
To understand our algorithm, Section \ref{section:convex_reformulation} describes the novel contributions that eventually lead to these new cuts, which are derived and applied in Section \ref{section:algorithm}.
Note that a more thorough algorithmic description is presented in Appendix \ref{section:appendix-algorithm}.


\section{Convex Reformulation}
\label{section:convex_reformulation}
Section \ref{subsec:cnlp} summarizes a convex method for determining the feasibility of a nonlinear network with fixed resistances.
This method could be exploited within a bilevel programming formulation for optimal design, whereby resistance selections obtained in the outer level must satisfy constraints on the corresponding solution of \eqref{eqn:p} (i.e., the inner level).
Indeed, such a bilevel method was developed by \citet{zhang1996bilevel} for the design of pipe networks.
Interestingly, however, no study has fully examined nor exploited the relationship between \eqref{eqn:p} and its (strong) dual, which we show leads to exact convex reformulations of the original feasibility and design problems.
This section describes two of our novel contributions, both of which serve as foundations for the remainder of this paper.
The first subsection derives an \emph{exact} convex reformulation of \eqref{eqn:nlp}.
The second subsection then extends this result to derive an \emph{exact} MICP reformulation of \eqref{eqn:minlp}.

\subsection{Convex Reformulation of (NLP(r))}
\label{subsection:nlp}
This section extends Section \ref{subsec:cnlp} to derive an exact convex reformulation of \eqref{eqn:nlp} based on optimality conditions for \eqref{eqn:p} and the addition of physical bounds.
Note that the objective function of \eqref{eqn:p} is convex and all constraints are affine.
The linearity constraint qualification thus implies the existence of a strong dual.
This dual problem is thoroughly derived in Appendix \ref{section:appendix-dual-derivation} using conventional Lagrangian duality theory.
It is expressed here as
\begin{equation}\tag{D$(r)$}\begin{aligned}
    & \underset{\Delta h^{\pm} \geq 0}{\text{maximize}}
    & & \frac{-\alpha}{1 + \alpha} \sum_{a \in \mathcal{A}} \frac{1}{\sqrt[\alpha]{L_{a} r_{a}}} \left[(\Delta h^{+}_{a})^{1 + \frac{1}{\alpha}} + (\Delta h^{-}_{a})^{1 + \frac{1}{\alpha}}\right] - \sum_{i \in \mathcal{J}} h_{i} d_{i} \\
    & \text{subject to}
    & & \Delta h_{a}^{+} - \Delta h_{a}^{-} = h_{i}^{s} - h_{j}, ~ \forall a := (i, j) \in \mathcal{A} : i \in \mathcal{S} \\
    & & & \Delta h_{a}^{+} - \Delta h_{a}^{-} = h_{i} - h_{j}, ~ \forall a := (i, j) \in \mathcal{A} : i \in \mathcal{J}.
\end{aligned}\label{eqn:d}\end{equation}

\begin{theorem}
Let $f_{P}(q)$ and $f_{D}(h)$ denote the objective functions of \eqref{eqn:p} and \eqref{eqn:d}, respectively, and $q_{a} := q_{a}^{+} - q_{a}^{-}, ~ \forall a \in \mathcal{A}$.
The following \emph{convex} problem is equivalent to \eqref{eqn:nlp}:
\begin{equation}\tag{CP$(r)$}\begin{aligned}
    & f_{P}(q) - f_{D}(h) \leq 0 \\
    & \sum_{\mathclap{a \in \delta^{-}_{i}}} q_{a} - \sum_{\mathclap{a \in \delta^{+}_{i}}} q_{a} = d_{i}, \; \forall i \in \mathcal{J} \\
    & \Delta h_{a}^{+} - \Delta h_{a}^{-} = h_{i}^{s} - h_{j}, ~ \forall a := (i, j) \in \mathcal{A} : i \in \mathcal{S} \\
    & \Delta h_{a}^{+} - \Delta h_{a}^{-} = h_{i} - h_{j}, ~ \forall a := (i, j) \in \mathcal{A} : i \in \mathcal{J} \\
    & \underline{h}_{i} \leq h_{i} \leq \overline{h}_{i}, ~ \forall i \in \mathcal{J}, ~ 0 \leq q_{a}^{\pm} \leq \overline{q}_{a}^{\pm}, ~ \Delta h_{a}^{\pm} \geq 0, ~ \forall a \in \mathcal{A}.
\end{aligned}\label{eqn:cp-exact}\end{equation}\label{thm:pd}\end{theorem}
\proof{Proof of Theorem \ref{thm:pd}}
By weak duality, it follows that $f_{P}(q) \geq f_{D}(h)$ for any feasible solutions to \eqref{eqn:p} and \eqref{eqn:d}, with equality holding for optimal solutions by strong duality.
As a result, optimality is equivalently imposed by combining constraints of \eqref{eqn:p} and \eqref{eqn:d} with the \emph{convex} constraint $f_{P}(q) - f_{D}(h) \leq 0$.
Next, \cite{raghunathan2013global} shows that head loss Constraints \eqref{eqn:feas-head-loss} of \eqref{eqn:nlp} are equivalent to a portion of the first order optimality conditions for \eqref{eqn:p}.
Moreover, Section 5.5.3 of \citet{boyd2004convex} ensures that, since \eqref{eqn:p} is a convex differentiable problem with a strong dual, the dual multipliers $h$ appearing in its optimality conditions are optimal solutions of the dual problem \eqref{eqn:d}.
That is, any $(q, h)$ that satisfies the constraints of \eqref{eqn:p} and \eqref{eqn:d} and strong duality will also satisfy flow conservation and head loss Constraints \eqref{eqn:feas-flow-conservation} and \eqref{eqn:feas-head-loss}.
Finally, appending bound constraints on $q$ and $h$ ensures equivalence to \eqref{eqn:nlp}.
\Halmos
\endproof

A physical interpretation of the convex strong duality constraint is discussed in Appendix \ref{section:appendix-strong-duality}.
Summarily, the constraint implies the conservation of power, with an inequality replacing the traditional equality.
Nonetheless, it is known via the strong duality argument above that any solution to \eqref{eqn:cp-exact} will indeed satisfy this constraint with equality.

\subsection{Mixed-integer Convex Reformulation of (MINLP)}
\label{subsection:micp}
To reformulate \eqref{eqn:minlp} using \eqref{eqn:cp-exact}, we first introduce continuous flow variables $q_{ap}^{\pm}$ as in Equations \eqref{eqn:flow-equations}, and binary resistance selection variables $x_{ap}$, subject to bound Constraints \eqref{eqn:micpr-flow-bounds}.
Next, we introduce continuous variables $\Delta h_{ap}^{\pm}$ and derive their bounds from Constraints \eqref{eqn:feas-head-bounds}, i.e.,
\begin{equation}
    0 \leq \Delta h^{\pm}_{ap} \leq \Delta \overline{h}^{\pm}_{ap} x_{ap}, ~ x_{ap} \in \mathbb{B}, ~ \forall a \in \mathcal{A}, ~ \forall p \in \mathcal{R}_{a} \label{eqn:micp-head-difference-bounds}.
\end{equation}
The constraints involving head differences in \eqref{eqn:cp-exact} are next rewritten as
\begin{subequations}
\begin{align}
    & \sum_{\mathclap{p \in \mathcal{R}_{a}}} \left(\Delta h_{ap}^{+} - \Delta h_{ap}^{-}\right) = h_{i} - h_{j}, ~ \forall a := (i, j) \in \mathcal{A} : i \in \mathcal{J} \\
    & \sum_{\mathclap{p \in \mathcal{R}_{a}}} \left(\Delta h_{ap}^{+} - \Delta h_{ap}^{-}\right) = h_{i}^{s} - h_{j}, ~ \forall a := (i, j) \in \mathcal{A} : i \in \mathcal{S}.
\end{align}\label{eqn:micp-head-equality}\end{subequations}
Finally, the strong duality constraint appearing in \eqref{eqn:cp-exact} is expanded over all $p \in \mathcal{R}_{a}$ as
\begin{equation}
\begin{gathered}
	\frac{1}{1 + \alpha} \sum_{\mathclap{a \in \mathcal{A}}} L_{a} \sum_{\mathclap{p \in \mathcal{R}_{a}}} p \left[(q_{ap}^{+})^{1 + \alpha} + (q_{ap}^{-})^{1 + \alpha}\right] - \sum_{i \in \mathcal{S}} h_{i}^{s} \sum_{a \in \delta^{+}_{i}} q_{a} \\
	+ ~ \frac{\alpha}{1 + \alpha} \sum_{a \in \mathcal{A}} \sum_{p \in \mathcal{R}_{a}} \frac{1}{\sqrt[\alpha]{L_{a} p}} \left[(\Delta h^{+}_{ap})^{1 + \frac{1}{\alpha}} + (\Delta h^{-}_{ap})^{1 + \frac{1}{\alpha}}\right] + \sum_{i \in \mathcal{J}} h_{i} d_{i} \leq 0.
\end{gathered}\label{eqn:micp-strong-duality}\end{equation}

These expansions of \eqref{eqn:cp-exact} give rise to the \emph{exact} mixed-integer convex reformulation of \eqref{eqn:minlp}, where Constraints \eqref{eqn:feas-flow-conservation} employ the definitions of $q_{a}$ described in Equations \eqref{eqn:flow-equations}.
Namely,
\begin{equation}\tag{MICP-E}\begin{aligned}
    & \text{minimize}
	& & \textnormal{Objective function:} ~ \eta(x) ~ \textnormal{of Equation} ~ \eqref{eqn:ne-obj} \\
    & \text{subject to}
	& & \textnormal{Physical bounds: Constraints} ~ \eqref{eqn:feas-head-bounds}, \eqref{eqn:micpr-flow-bounds}, \eqref{eqn:micp-head-difference-bounds} \\
	& & & \textnormal{Flow conservation: Constraints} ~ \eqref{eqn:feas-flow-conservation} \\
	& & & \textnormal{Resistance selection: Constraints} ~ \eqref{eqn:ne-resistance-selection} \\
	& & & \textnormal{Head difference equalities: Constraints} ~ \eqref{eqn:micp-head-equality} \\
	& & & \textnormal{Strong duality: Constraint} ~ \eqref{eqn:micp-strong-duality}.
\end{aligned}\label{eqn:micp}\end{equation}

\begin{theorem}
\label{eqn:thm-equiv}
A design $\hat{x}$ is feasible for \eqref{eqn:micp} if and only if it is feasible for \eqref{eqn:minlp}.
\end{theorem}
\proof{Proof of Theorem \ref{eqn:thm-equiv}}
For any binary $\hat{x}$ satisfying Constraints \eqref{eqn:ne-resistance-selection}, \eqref{eqn:minlp} reduces to \eqref{eqn:nlp} and \eqref{eqn:micp} to \eqref{eqn:cp-exact}, with $r$ given by Equations \eqref{eqn:x-relation}.
For any $r$, \eqref{eqn:nlp} and \eqref{eqn:cp-exact} are equivalent by Theorem \ref{thm:pd}.
Thus, the sets of feasible $\hat{x}$ for \eqref{eqn:minlp} and \eqref{eqn:micp} are equal.
\Halmos
\endproof

Constraint \eqref{eqn:micp-strong-duality} of \eqref{eqn:micp} can be viewed as a convex embedding of Constraints \eqref{eqn:minlp-head-loss}.
Although convexity is desirable, Constraint \eqref{eqn:micp-strong-duality} is highly aggregated.
In this sense, the disaggregated nonconvex Constraints \eqref{eqn:minlp-head-loss} of \eqref{eqn:minlp} or convex Constraints \eqref{eqn:micpr-head-loss} of \eqref{eqn:micpr} may be more numerically useful.
Section \ref{section:algorithm} uses this observation to construct a relaxation-based global algorithm based on \eqref{eqn:micp} and \eqref{eqn:micpr} that outer-approximates Constraint \eqref{eqn:micp-strong-duality} \emph{and} Constraints \eqref{eqn:micpr-head-loss}.


\section{Novel Cutting Planes and Algorithmic Enhancements}
\label{section:algorithm}
As discussed in Section \ref{subsection:micp}, \eqref{eqn:micp} is an exact reformulation of the design problem.
As such, it can be posed directly to an MICP solver (e.g., \textsc{BONMIN}, \citealp{bonami2008algorithmic}) to obtain a globally optimal solution.
In practice, however, direct methods converge slowly, as modern MICP solvers do not efficiently handle nonquadratic nonlinear relationships (e.g., Constraint \eqref{eqn:micp-strong-duality}).
On the other hand, modern MIP solvers are highly efficient but require conscientious linearizations of \eqref{eqn:micp}.
This section pursues the latter, following a structure similar to Section \ref{section:relaxation_based_reformulation}.
Specifically, it introduces novel outer approximation cuts while developing a MIP reformulation and relaxation of \eqref{eqn:micp}, then summarizes this paper's extensions to the algorithm of \cite{raghunathan2013global}.

\paragraph{Flow Direction-based Inequalities}
Although \eqref{eqn:micp} does not require flow direction variables $y_{a}$, we nonetheless incorporate $y_{a}$ to strengthen inequalities throughout our linearized reformulation.
Similar to Constraints \eqref{eqn:micpr-direction-bounds}, flows and head differences are bounded via
\begin{subequations}\begin{gather}
	0 \leq q_{ap}^{+} \leq \overline{q}_{ap}^{+} y_{a}, ~ 0 \leq q_{ap}^{-} \leq \overline{q}_{ap}^{-} (1 - y_{a}), ~ y_{a} \in \mathbb{B}, ~ \forall a \in \mathcal{A}, ~ \forall p \in \mathcal{R}_{a} \\
	0 \leq \Delta h_{ap}^{+} \leq \Delta \overline{h}_{ap}^{+} y_{a}, ~ 0 \leq \Delta h_{ap}^{-} \leq \Delta \overline{h}_{ap}^{-} (1 - y_{a}), ~ y_{a} \in \mathbb{B}, ~ \forall a \in \mathcal{A}, ~ \forall p \in \mathcal{R}_{a}.
\end{gather}\label{eqn:micpe-direction-bounds}\end{subequations}
We also employ valid inequalities to exploit a priori knowledge concerning flow directionality throughout the network.
The inequalities proposed here are similar to those described by \cite{borraz2016convex} in the context of natural gas network expansion planning.
The first are
\begin{equation}
    \sum_{a \in \delta_{i}^{+}} y_{a} \geq 1, ~ \forall i \in \mathcal{S} \label{eqn:source-flow},
\end{equation}
which model that at least one pipe must send water \emph{away} from each source.
The next are
\begin{equation}
    \sum_{a \in \delta_{i}^{-}} y_{a} + \sum_{a \in \delta_{i}^{+}} (1 - y_{a}) \geq 1, ~ i \in \mathcal{J} : d_{i} > 0 \label{eqn:demand-flow},
\end{equation}
which model that at least one pipe must provide water \emph{to} each demand node.
Finally,
\begin{subequations}\begin{align}
	& \sum_{a \in \delta_{i}^{-}} y_{a} - \sum_{a \in \delta_{i}^{+}} y_{a} = 0, ~ i \in \mathcal{J} : \left(d_{i} = 0\right) \land \left(\deg^{+}_{i} = \deg^{-}_{i} = 1\right) \label{eqn:deg-2-flow-1} \\
	& \sum_{a \in \delta_{i}^{-}} y_{a} + \sum_{a \in \delta_{i}^{+}} y_{a} = 1, ~ i \in \mathcal{J} : \left(d_{i} = 0\right) \land \left(\deg^{\pm}_{i} = 2\right) \land \left(\deg^{\mp}_{i} = 0\right) \label{eqn:deg-2-flow-2}
\end{align}\label{eqn:deg-2-flow}\end{subequations}
are constraints that model flow directionality at junctions with zero demand and degree two, with the implication that the direction of incoming flow must be equal to the direction of outgoing flow.

\paragraph{Head Loss Outer Approximation Cutting Planes}
Although also not required by \eqref{eqn:micp}, head loss relationships are used to further strengthen linearization-based reformulations.
Similar to Constraints \eqref{eqn:imip-head-loss-aggregated-4} in \eqref{eqn:milpr}, we express these head loss outer approximations as
\begin{subequations}\begin{align}
    \tilde{\tau}_{ar_{a}}^{+} y_{a} + \alpha \sum_{\mathclap{p \in \mathcal{R}_{a}}} p \left(\tilde{q}_{ap}^{+}\right)^{\alpha - 1} q_{ap}^{+} &\leq \sum_{\mathclap{p \in \mathcal{R}_{a}}} \frac{\Delta h_{ap}^{+}}{L_{a}}, \, \forall a \in \mathcal{A}, \, \forall r_{a} \in \mathcal{R}_{a}, \, \forall \tilde{q}_{ar_{a}}^{+} \in \mathcal{Q}^{+}_{ar_{a}} \\
    \tilde{\tau}_{ar_{a}}^{-} \left(1 - y_{a}\right) + \alpha \sum_{\mathclap{p \in \mathcal{R}_{a}}} p \left(\tilde{q}_{ap}^{-}\right)^{\alpha - 1} q_{ap}^{-} &\leq \sum_{\mathclap{p \in \mathcal{R}_{a}}} \frac{\Delta h_{ap}^{-}}{L_{a}}, \, \forall a \in \mathcal{A}, \, \forall r_{a} \in \mathcal{R}_{a}, \, \forall \tilde{q}_{ar_{a}}^{-} \in \mathcal{Q}^{-}_{ar_{a}}.
\end{align}\label{eqn:mip-head-loss}\end{subequations}
Furthermore, similar to Constraints \eqref{eqn:micpr-head-loss}, we linearly upper-bound each head difference with
\begin{equation}
    \Delta h_{ap}^{\pm} \leq L_{a} p \left(\overline{q}_{ap}^{\pm}\right)^{\alpha - 1} q_{ap}^{\pm}, ~ \forall a \in \mathcal{A}, ~ \forall p \in \mathcal{R}_{a} \label{eqn:mip-head-difference-ub}.
\end{equation}

\paragraph{Strong Duality Cutting Planes}
A primary contribution of this paper is the strong duality Constraint \eqref{eqn:micp-strong-duality} of \eqref{eqn:micp}.
To linearize this constraint, we define  variables $q_{a}^{\textrm{NL}}, \Delta h_{a}^{\textrm{NL}} \geq 0$ to approximate the sums of nonlinear terms involving $q_{ap}^{\pm}$ and $\Delta h_{ap}^{\pm}$ in the original constraint, i.e.,
\begin{subequations}
\begin{align}
    q_{a}^{\textrm{NL}} &= \frac{1}{1 + \alpha} \sum_{p \in \mathcal{R}_{a}} p \left[\left(q_{ap}^{+}\right)^{1 + \alpha} + \left(q_{ap}^{-}\right)^{1 + \alpha} \right], \forall a \in \mathcal{A} \\
    \Delta h_{a}^{\textrm{NL}} &= \frac{\alpha}{1 + \alpha} \sum_{p \in \mathcal{R}_{a}} \frac{1}{\sqrt[\alpha]{p}} \left[(\Delta h^{+}_{ap})^{1 + \frac{1}{\alpha}} + (\Delta h^{-}_{ap})^{1 + \frac{1}{\alpha}}\right], ~ \forall a \in \mathcal{A}.
\end{align}
\label{eqn:nl-relationships}
\end{subequations}
This permits a purely linear rewriting of the original strong duality constraint, namely
\begin{equation}
	\sum_{\mathclap{a \in \mathcal{A}}} L_{a} q_{a}^{\textrm{NL}} - \sum_{i \in \mathcal{S}} h_{i}^{s} \sum_{a \in \delta^{+}_{i}} \sum_{p \in \mathcal{R}_{a}} (q_{ap}^{+} - q_{ap}^{-}) + \sum_{a \in \mathcal{A}} \frac{\Delta h_{a}^{\textrm{NL}}}{\sqrt[\alpha]{L_{a}}} + \sum_{i \in \mathcal{J}} h_{i} d_{i} \leq 0.
\label{eqn:mip-strong-duality}\end{equation}

Observe that the right-hand sides of Equations \eqref{eqn:nl-relationships} are convex.
Similar to the head loss outer approximations derived in Section \ref{subsec:mip-relaxation}, we follow the notion of \emph{equal intercept} linear outer approximations to compose cuts similar to Constraints \eqref{eqn:mip-head-loss} for $q_{a}^{\textrm{NL}}$ and $\Delta h_{a}^{\textrm{NL}}$.
For $q_{a}^{\textrm{NL}}$, let the intercept be determined by the linear approximation to the term corresponding to $r_{a} \in \mathcal{R}_{a}$ at point $\tilde{q}_{ar_{a}}^{\pm}$.
For each remaining $p \in \mathcal{R}_{a}$, the same intercept is achieved by the linear approximation at point $\tilde{q}_{ap}^{\pm} := \tilde{q}_{ar_{a}}^{\pm}(r_{a}  / p)^{1 / (1 + \alpha)}$.
With these values in mind, the outer approximations are then
\begin{subequations}\begin{align}
    \tilde{\zeta}_{ar_{a}}^{+} y_{a} + \sum_{\mathclap{p \in \mathcal{R}_{a}}} p \left(\tilde{q}_{ap}^{+}\right)^{\alpha} q_{ap}^{+} &\leq q_{a}^{\textrm{NL}}, \, \forall a \in \mathcal{A}, \, \forall r_{a} \in \mathcal{R}_{a}, \, \forall \tilde{q}_{ar_{a}}^{+} \in \mathcal{Q}^{\textrm{NL}^{+}}_{ar_{a}} \\
    \tilde{\zeta}_{ar_{a}}^{-} \left(1 - y_{a}\right) + \sum_{\mathclap{p \in \mathcal{R}_{a}}} p \left(\tilde{q}_{ap}^{-}\right)^{\alpha} q_{ap}^{-} &\leq q_{a}^{\textrm{NL}}, \, \forall a \in \mathcal{A}, \, \forall r_{a} \in \mathcal{R}_{a}, \, \forall \tilde{q}_{ar_{a}}^{-} \in \mathcal{Q}^{\textrm{NL}^{-}}_{ar_{a}},
\end{align}\label{eqn:mip-q-nl-cuts}\end{subequations}
where $\mathcal{Q}^{\textrm{NL} \pm}_{ar_{a}} = [0, \overline{q}_{ar_{a}}^{\pm}]$, and $\tilde{\zeta}_{ar_{a}}^{\pm} := \left(\frac{1}{1 + \alpha} - 1\right) r_{a} \left(\tilde{q}_{ar_{a}}^{\pm}\right)^{1 + \alpha}$ is defined for notational brevity.

An analogous procedure is followed for $\Delta h_{a}^{\textrm{NL}}$.
Similar to Constraints \eqref{eqn:mip-q-nl-cuts}, assuming $\Delta \tilde{h}_{ap}^{\pm}$ are points with coinciding outer approximation intercepts, the aggregated outer approximations are
\begin{subequations}\begin{align}
    -\xi_{ar_{a}}^{\pm} y_{a} + \sum_{\mathclap{p \in \mathcal{R}_{a}}} \left(\sqrt[\alpha]{\frac{\Delta \tilde{h}_{ap}^{+}}{p}}\right) \Delta h_{ap}^{+}
    &\leq \Delta h_{a}^{\textrm{NL}}, \, \forall a \in \mathcal{A}, \, \forall r_{a} \in \mathcal{R}_{a}, \, \forall \Delta \tilde{h}_{ar_{a}}^{+} \in \mathcal{H}^{\textrm{NL}^{+}}_{ar_{a}} \\
    -\xi_{ar_{a}}^{\pm} (1 - y_{a}) + \sum_{\mathclap{p \in \mathcal{R}_{a}}} \left(\sqrt[\alpha]{\frac{\Delta \tilde{h}_{ap}^{-}}{p}}\right) \Delta h_{ap}^{-}
    &\leq \Delta h_{a}^{\textrm{NL}}, \, \forall a \in \mathcal{A}, \, \forall r_{a} \in \mathcal{R}_{a}, \, \forall \Delta \tilde{h}_{ar_{a}}^{-} \in \mathcal{H}^{\textrm{NL}^{-}}_{ar_{a}},
\end{align}\label{eqn:mip-dh-nl-cuts}\end{subequations}
where $\mathcal{H}^{\textrm{NL}^{\pm}}_{ar_{a}} = [0, \Delta \overline{h}_{ar_{a}}^{\pm}]$, and $\xi_{ar_{a}}^{\pm} := \frac{\left(\Delta \tilde{h}_{ar_{a}}^{\pm}\right)^{1 + \frac{1}{\alpha}}}{\left(1 + \alpha\right)\sqrt[\alpha]{r_{a}}}$ is similarly defined for notational brevity.

\paragraph{Mixed-integer Linear Reformulation}
With the variables and constraints previously described, a MIP reformulation of \eqref{eqn:micp} may be written in a manner similar to \eqref{eqn:milpr}, that is,
\begin{equation}\tag{MIP-E}\begin{aligned}
    & \text{minimize}
	& & \textnormal{Objective function:} ~ \eta(x) ~ \textnormal{of Equation} ~ \eqref{eqn:ne-obj} \\
    & \text{subject to}
	& & \textnormal{Physical bounds: Constraints} ~ \eqref{eqn:feas-head-bounds}, \eqref{eqn:micpr-flow-bounds}, \eqref{eqn:micp-head-difference-bounds} \\
	& & & \textnormal{Flow conservation: Constraints} ~ \eqref{eqn:feas-flow-conservation} \\
	& & & \textnormal{Resistance selection: Constraints} ~ \eqref{eqn:ne-resistance-selection} \\
	& & & \textnormal{Head difference relationships: Constraints} ~ \eqref{eqn:micp-head-equality}, \eqref{eqn:mip-head-loss}, \eqref{eqn:mip-head-difference-ub} \\
    & & & \textnormal{Direction-related inequalities: Constraints} ~ \eqref{eqn:micpe-direction-bounds}\text{--}\eqref{eqn:deg-2-flow} \\
    & & & \textnormal{Strong duality: Constraints} ~ \eqref{eqn:mip-q-nl-cuts}, \eqref{eqn:mip-dh-nl-cuts} \\
    & & & \textnormal{Feasibility cutting planes: Constraints} ~ \eqref{eqn:no-good-cut}.
\end{aligned}\label{eqn:mip}\end{equation}
Similarly to \eqref{eqn:milpr}, algorithmically, the sets $\mathcal{Q}^{\pm}$, $\mathcal{Q}^{\textrm{NL} \pm}$, $\mathcal{H}^{\textrm{NL} \pm}$, and $\bar{\mathcal{X}}$ of Constraints \eqref{eqn:mip-head-loss}, \eqref{eqn:mip-q-nl-cuts}, \eqref{eqn:mip-dh-nl-cuts}, and \eqref{eqn:no-good-cut} are instead replaced by the initially empty finite sets $\tilde{\mathcal{Q}}^{\pm}$, $\tilde{\mathcal{Q}}^{\textrm{NL} \pm}$, $\tilde{\mathcal{H}}^{\textrm{NL} \pm}$, and $\tilde{\mathcal{X}}$, respectively.
These changes give rise to the further (finite) linear relaxation, named (MIP-ER).

\paragraph{Algorithmic Enhancements}
Our algorithm, which is omitted here for brevity but exploits (MIP-ER), is similar to that of \cite{raghunathan2013global} but differs in a few important respects.
Primarily, Raghunathan only applies outer approximations similar to Constraints \eqref{eqn:mip-head-loss}, which correspond to convexified head loss relationships.
Our algorithm extends this by adding outer approximations of terms appearing in the strong duality Constraint \eqref{eqn:micp-strong-duality}.
A more thorough comparison and presentation of the two algorithms is presented in Appendix \ref{section:appendix-algorithm}.


\section{Computational Experiments}
\label{section:computational_experiments}
This section compares the convergence of our new algorithm and an algorithm based on \citet{raghunathan2013global}.
Both were implemented in the \textsc{Julia} programming language using \textsc{JuMP}, version 0.20 \citep{dunning2017jump}, and version 0.1 of \textsc{WaterModels}, an open-source \textsc{Julia} package for water distribution network optimization \citep{watermodels}.
Section \ref{subsec:experimental_setup} describes the instances, computational resources, and parameters used in the experiments; Section \ref{subsec:comparison} compares the efficacy of the two algorithms by examining convergence; and Section \ref{subsec:scaling} compares their performance on an extended set of instances obtained by scaling demand across all junctions.

\subsection{Experimental Setup}
\label{subsec:experimental_setup}
The numerical experiments consider instances of varying sizes that appear in the water network design literature and are summarized in Table \ref{table:problem_summary} \citep{d2015mathematical,mala2018lost}.
All use the Hazen-Williams head loss relationship originally defined by Equation \eqref{eqn:hazen-williams}.
The set of \emph{diameters} in each problem gives a set of \emph{resistances}, each being proportional to $D_{ap}^{-4.8704}$.

{\renewcommand{\arraystretch}{0.5}
\begin{table}[t]
    \begin{center}
        \begin{tabular}{|c|c|c|c|c|c|}
            \hline
            Network & \# Nodes & \# Arcs & \# Resistances & \# Binary Variables & \# Designs \\ \hline
            \texttt{shamir}        & 8   & 8   & 14 & 120  & $1.48 \times 10^{9}$   \\ \hline
            \texttt{blacksburg}    & 32  & 23  & 14 & 369  & $2.30 \times 10^{26}$  \\ \hline
            \texttt{hanoi}         & 33  & 34  & 6  & 238  & $2.87 \times 10^{26}$  \\ \hline
            \texttt{foss\_poly\_0} & 38  & 58  & 7  & 464  & $1.04 \times 10^{49}$  \\ \hline
            \texttt{foss\_iron}    & 38  & 58  & 13 & 812  & $4.06 \times 10^{64}$  \\ \hline
            \texttt{foss\_poly\_1} & 38  & 58  & 22 & 1334 & $7.25 \times 10^{77}$  \\ \hline
            \texttt{pescara}       & 74  & 99  & 13 & 1386 & $1.91 \times 10^{110}$ \\ \hline
            \texttt{modena}        & 276 & 317 & 13 & 4438 & $1.32 \times 10^{353}$ \\ \hline
        \end{tabular}
    \end{center}
    \caption{Summary of benchmark instances from the literature and their combinatorial sizes. Here, ``\# Arcs'' is the number of arcs with $\lvert \mathcal{R}_{a} \rvert \neq 1$; ``\# Resistances'' is $\lvert \mathcal{R}_{a} \rvert$ per \emph{variable} pipe; ``\# Binary Variables'' is $\lvert x \rvert + \lvert y \rvert$ for MIP formulations; and ``\# Designs'' is the total number of unique designs satisfying Constraints \eqref{eqn:ne-resistance-selection}.}
    \label{table:problem_summary}
\end{table}}

The instances of Table \ref{table:problem_summary} are divided into two classes: \emph{moderate} instances, comprising \texttt{shamir}, \texttt{blacksburg}, \texttt{hanoi}, \texttt{foss\_poly\_0}, and \texttt{foss\_iron}; and \emph{large} instances, comprising \texttt{foss\_poly\_1}, \texttt{pescara}, and \texttt{modena}.
Generally, moderate instances are solvable to optimality with both algorithms given a sufficient amount of time (i.e., seconds to hours), while large instances cannot be solved, even given substantial time (i.e., days).
Each experiment began each algorithm with equivalent data, initial feasible solutions, and outer approximation points.
Parameters of the two algorithms were chosen to coincide with those used by Raghunathan and are detailed in Appendix \ref{appendix:nodecuts}.

Experiments were performed on Los Alamos National Laboratory's Darwin computing cluster.
Each was executed on a node containing two Intel Xeon E5-2695 v4 processors, each with 18 cores @2.10 GHz, and 125 GB of memory.
Excluding the small amount of time required by Raghunathan's heuristic procedure to obtain an initial feasible solution (seconds), each experiment was provided a wall-clock time of $171{,}900$ seconds (approximately two days).
For solutions of the MIPs, \textsc{Gurobi} 9.0.3 was used with \texttt{Cuts=0}, which disables all of \textsc{Gurobi}'s internal cutting plane methods.
For moderate instances, \texttt{Heuristics=0.0} was used, which disables \textsc{Gurobi}'s internal heuristics.
For large instances, \texttt{MIPFocus=1} was used, which places a focus on finding feasible solutions quickly.

For convex subproblems, e.g., the solution of \eqref{eqn:p}, \textsc{Ipopt} version 3.13 was used \citep{wachter2006implementation}.
As per \cite{tasseff2019exploring}, since these problems are small, the linear solver \textsc{MA57} was employed.
The settings \texttt{warm\_start\_init\_point="yes"} and \texttt{nlp\_scaling\_method="none"} were also used.
Heuristically, these parameters computed solutions to \eqref{eqn:p} most efficiently.

\subsection{Comparison of Algorithms on Standard Benchmark Instances}
\label{subsec:comparison}
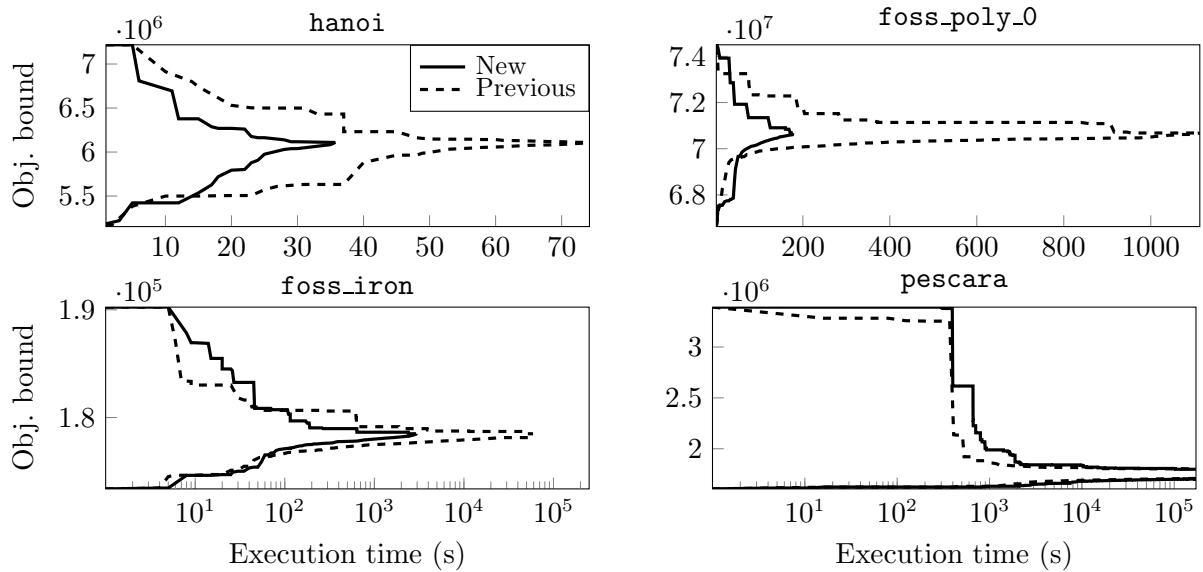
\begin{figure}[!ht]
    \centering
    \begin{subfigure}[b]{0.49\textwidth}
        \centering
        \begin{tikzpicture}[]
	\centering
	\begin{axis}[legend cell align=left,enlargelimits=false,xtick pos=left,
                 ytick pos=left,height=4.0cm,ylabel=Obj. bound,
                 width=0.99\textwidth,legend style={at={(1.00, 1.00)},anchor=north east},
				 legend columns=1,title=\texttt{hanoi},title style={yshift=-1.0ex},
				 legend style={font=\small,column sep=1.5pt, row sep=-4.0pt}]
		\pgfplotstableread[col sep = comma]{computational_experiments/data/standard/NETWORK.hanoi-SIAM.false.csv}{\tasseff};
		\addplot[very thick, opacity=1.0] table [x = time_elapsed, y = lower_bound]{\tasseff};
		\addplot[very thick, opacity=1.0, forget plot] table [x = time_elapsed, y = upper_bound]{\tasseff};
		\pgfplotstableread[col sep = comma]{computational_experiments/data/standard/NETWORK.hanoi-SIAM.true.csv}{\siam};
		\addplot[very thick, dashed] table [x = time_elapsed, y = lower_bound]{\siam};
		\addplot[very thick, dashed, forget plot] table [x = time_elapsed, y = upper_bound]{\siam};
		\legend{New,Previous};
	\end{axis}
\end{tikzpicture}
    \end{subfigure}
    \hfill
    \begin{subfigure}[b]{0.49\textwidth}
        \centering
        \begin{tikzpicture}
	\centering
	\begin{axis}[legend cell align=left,enlargelimits=false,xtick pos=left,
                 ymax=7.45e7,ytick pos=left,height=4.0cm,width=0.99\textwidth,
                 /pgf/number format/.cd,1000 sep={},
                 title=\texttt{foss\_poly\_0},
                 title style={yshift=-1.0ex}]
		\pgfplotstableread[col sep = comma]{computational_experiments/data/standard/NETWORK.foss_poly_0-SIAM.false.csv}{\tasseff};
		\addplot[very thick, opacity=1.0] table [x = time_elapsed, y = lower_bound]{\tasseff};
		\addplot[very thick, opacity=1.0, forget plot] table [x = time_elapsed, y = upper_bound]{\tasseff};
		\pgfplotstableread[col sep = comma]{computational_experiments/data/standard/NETWORK.foss_poly_0-SIAM.true.csv}{\siam};
		\addplot[very thick, dashed, opacity=1.0] table [x = time_elapsed, y = lower_bound]{\siam};
		\addplot[very thick, dashed, opacity=1.0, forget plot] table [x = time_elapsed, y = upper_bound]{\siam};
	\end{axis}
\end{tikzpicture}
    \end{subfigure} \\
    \begin{subfigure}[b]{0.49\textwidth}
        \centering
        \begin{tikzpicture}[]
	\centering
	\begin{semilogxaxis}[legend cell align=left,enlargelimits=false,xtick pos=left,xmax=2.5e5,
                 ytick pos=left,height=4.0cm,ylabel=Obj. bound,
                 xlabel=Execution time (s),
                 ytick={1.8e5,1.9e5},xtick={1.0e1,1.0e2,1.0e3,1.0e4,1.0e5},
                 extra x ticks={2,3,4,5,6,7,8,9,2.0e5},
                 extra x tick style={tickwidth=\pgfkeysvalueof{/pgfplots/minor tick length},xticklabels={}},
                 width=0.99\textwidth,title=\texttt{foss\_iron},title style={yshift=-1.0ex}]
		\pgfplotstableread[col sep = comma]{computational_experiments/data/standard/NETWORK.foss_iron-SIAM.false.csv}{\tasseff};
		\addplot[very thick, opacity=1.0] table [x = time_elapsed, y = lower_bound]{\tasseff};
		\addplot[very thick, opacity=1.0, forget plot] table [x = time_elapsed, y = upper_bound]{\tasseff};
		\pgfplotstableread[col sep = comma]{computational_experiments/data/standard/NETWORK.foss_iron-SIAM.true.csv}{\siam};
		\addplot[very thick, dashed, opacity=1.0] table [x = time_elapsed, y = lower_bound]{\siam};
		\addplot[very thick, dashed, opacity=1.0, forget plot] table [x = time_elapsed, y = upper_bound]{\siam};
	\end{semilogxaxis}
\end{tikzpicture}
    \end{subfigure}
    \hfill
    \begin{subfigure}[b]{0.49\textwidth}
        \centering
        \begin{tikzpicture}[]
	\centering
	\begin{semilogxaxis}[legend cell align=left,enlargelimits=false,xtick pos=left,ymax=3392146.46,
                 ytick pos=left,height=4.0cm,xlabel=Execution time (s),
	             width=0.99\textwidth,legend style={at={(1.00, 1.00)},anchor=north east},
                 xtick={1.0e1,1.0e2,1.0e3,1.0e4,1.0e5},
                 extra x ticks={2,3,4,5,6,7,8,9,2.0e5},
                 extra x tick style={tickwidth=\pgfkeysvalueof{/pgfplots/minor tick length},xticklabels={}},
				 legend columns=1,title=\texttt{pescara},title style={yshift=-1.0ex},
				 legend style={font=\small,column sep=1.5pt, row sep=-4.0pt}]
		\pgfplotstableread[col sep = comma]{computational_experiments/data/standard/NETWORK.pescara-SIAM.false.csv}{\tasseff};
		\addplot[very thick, opacity=1.0] table [x = time_elapsed, y = lower_bound]{\tasseff};
		\addplot[very thick, opacity=1.0, forget plot] table [x = time_elapsed, y = upper_bound]{\tasseff};
		\pgfplotstableread[col sep = comma]{computational_experiments/data/standard/NETWORK.pescara-SIAM.true.csv}{\siam};
		\addplot[very thick, dashed, opacity=1.0] table [x = time_elapsed, y = lower_bound]{\siam};
		\addplot[very thick, dashed, opacity=1.0, forget plot] table [x = time_elapsed, y = upper_bound]{\siam};
	\end{semilogxaxis}
\end{tikzpicture}
    \end{subfigure}
    \caption{Convergence of objective bounds on select instances, comparing the new algorithm with an algorithm similar to \cite{raghunathan2013global} (Previous). Note the use of linear and logarithmic abscissas.}
    \label{fig:standard-summary}
\end{figure}

{\renewcommand{\arraystretch}{0.5}
\begin{table}[t]
    \begin{center}
    \begin{tabular}{c|c|c|c|c|c|c|}
        \cline{2-7} & \multicolumn{3}{c|}{Previous Algorithm} & \multicolumn{3}{c|}{New Algorithm} \\
        \cline{1-7} \multicolumn{1}{|c|}{Problem} & Gap (\%) & Nodes Expl. & Time (s) & Gap (\%) & Nodes Expl. & Time (s) \\
        \cline{1-7} \multicolumn{1}{|c|}{\texttt{shamir}}        & $0.00$ & $12{,}098$ & $12.27$   & $0.00$          & $\bm{2{,}567}$ & $\bm{7.11}$ \\
        \cline{1-7} \multicolumn{1}{|c|}{\texttt{hanoi}}         & $0.00$ & $32{,}024$ & $74.21$  & $0.00$          & $\bm{24{,}765}$ & $\bm{35.74}$ \\
        \cline{1-7} \multicolumn{1}{|c|}{\texttt{blacksburg}}    & $0.00$ & $16{,}009$ & $\bm{14.05}$ & $0.00$          & $\bm{15{,}971}$ & $29.25$ \\
        \cline{1-7} \multicolumn{1}{|c|}{\texttt{foss\_poly\_0}} & $0.00$ & $144{,}226$ & $1{,}112.78$        & $0.00$          & $\bm{63{,}120}$ & $\mathbf{177.08}$ \\
        \cline{1-7} \multicolumn{1}{|c|}{\texttt{foss\_iron}}    & $0.00$ & $1{,}307{,}123$ & $59{,}088.98$             & $0.00$ & $\bm{282{,}202}$ & $\bm{2{,}923.24}$ \\
        \cline{1-7} \multicolumn{1}{|c|}{\texttt{foss\_poly\_1}} & $\bm{4.19}$ & $48{,}320{,}343$ & Limit             & $4.90$ & $21{,}858{,}635$ & Limit \\
        \cline{1-7} \multicolumn{1}{|c|}{\texttt{pescara}}       & $\bm{5.26}$ & $5{,}633{,}461$ & Limit             & $5.29$ & $2{,}010{,}998$ & Limit \\
        \cline{1-7} \multicolumn{1}{|c|}{\texttt{modena}}        & $\bm{33.77}$ & $329{,}614$ & Limit             & $41.65$ & $55{,}592$ & Limit \\
        \cline{1-7}
    \end{tabular}
    \end{center}
    \caption{Comparison of optimality gaps, nodes explored, and solution times for the new algorithm and one similar to \cite{raghunathan2013global} (Previous). Bold denotes better (smaller) times, nodes explored, and gaps.}
    \label{tab:comparison}
\end{table}}

{\renewcommand{\arraystretch}{0.5}
\begin{table}[!ht]
    \begin{center}
    \begin{tabular}{c|c|c|c|c|}
        \cline{2-5} & \multicolumn{2}{c|}{Solutions from the Literature} & \multicolumn{2}{c|}{Solutions from This Study} \\
        \cline{1-5} \multicolumn{1}{|c|}{Problem} & Lower Bnd. & Upper Bnd. & Lower Bnd. & Upper Bnd. \\
        \cline{1-5} \multicolumn{1}{|c|}{\texttt{shamir}}        & $\bm{419{,}000}^{1*}$ & $\bm{419{,}000}^{1*}$ & $\bm{419{,}000}^{*}$ & $\bm{419{,}000}^{*}$ \\
        \cline{1-5} \multicolumn{1}{|c|}{\texttt{hanoi}}         & $\bm{6{,}109{,}620.09}^{1*}$ & $\bm{6{,}109{,}620.09}^{1*}$ & $\bm{6{,}109{,}620.90}^{*}$ & $\bm{6{,}109{,}620.90}^{*}$ \\
        \cline{1-5} \multicolumn{1}{|c|}{\texttt{blacksburg}}    & $\bm{118{,}251.09}^{1*}$ & $\bm{118{,}251.09}^{1*}$ & $\bm{118{,}251.09}^{*}$ & $\bm{118{,}251.09}^{*}$ \\
        \cline{1-5} \multicolumn{1}{|c|}{\texttt{foss\_poly\_0}} & $70{,}063{,}161.90^{1}$ & $\bm{70{,}680{,}507.90}^{2}$ & $\color{palette2}{\bm{70{,}680{,}507.90}}^{*}$ & $\color{palette2}{\bm{70{,}680{,}507.90}}^{*}$ \\
        \cline{1-5} \multicolumn{1}{|c|}{\texttt{foss\_iron}}    & $177{,}512.42^{1}$ & $\bm{178{,}494.14}^{1}$ & $\color{palette2}{\bm{178{,}494.14}}^{*}$ & $\color{palette2}{\bm{178{,}494.14}}^{*}$ \\
        \cline{1-5} \multicolumn{1}{|c|}{\texttt{foss\_poly\_1}} & $26{,}240.84^{1}$ & $29{,}202.99^{2}$ & $\bm{27{,}269.65}$ & $\bm{28{,}462.34}$ \\
        \cline{1-5} \multicolumn{1}{|c|}{\texttt{pescara}}       & $1{,}700{,}517.06^{1}$ & $\approx \bm{1{,}790{,}000}^{3}$ & $\bm{1{,}708{,}090.52}$ & $1{,}798{,}252.52$ \\
        \cline{1-5} \multicolumn{1}{|c|}{\texttt{modena}}        & $\bm{2{,}206{,}914.89}^{1}$ & $\approx \bm{2{,}560{,}000}^{4}$ & $2{,}198{,}756.06$ & $3{,}319{,}652.71$ \\
        \cline{1-5}
    \end{tabular}
    \end{center}
    \caption{Comparison of best lower and upper bounds from the literature and this study, with bold denoting best bounds, asterisks denoting proven optimality, and blue denoting instances closed for the first time. References from the literature are labeled as $^{1}$\citep{raghunathan2013global}, $^{2}$\citep{bragalli2012optimal}, $^{3}$\citep{zheng2017adaptive}, and $^{4}$\citep{artina2012contribution}.}
    \label{tab:best-bounds}
\end{table}}

Figure \ref{fig:standard-summary} illustrates the lower and upper bound convergence of the two algorithms on a representative subset of instances (i.e., \texttt{hanoi}, \texttt{foss\_poly\_0}, \texttt{foss\_iron}, and \texttt{pescara}).
Here, both algorithms converge to global optimality on the three moderate instances, with the new algorithm displaying more favorable performance, i.e., lower and upper bounds converging more quickly.
However, both algorithms can solve \texttt{hanoi} in relatively short amounts of times, taking just over a minute to reach global optimality in the worst case.
For larger moderate instances (i.e., \texttt{foss\_poly\_0} and \texttt{foss\_iron}), the differences in convergence behavior are more dramatic.
In the case of \texttt{foss\_poly\_0}, the new algorithm converges nearly an order of magnitude more quickly.
For \texttt{foss\_iron}, the difference is further emphasized, with the new algorithm converging more than an order of magnitude faster.
\emph{To highlight this point, throughout the literature, neither of these two \texttt{foss} instances appear to have been solved to global optimality.
Here, both algorithms solve both problems, but ours does so around an order of magnitude faster.}
In large cases (e.g., \texttt{pescara}), the algorithm of \cite{raghunathan2013global} converges more quickly.
This could be for many reasons, the most likely being the more frequent addition of outer approximations in the new algorithm, creating larger BB linear subproblems.

Table \ref{tab:comparison} provides relevant convergence data for the instances.
This table further supports the trends exemplified by Figure \ref{fig:standard-summary}.
For moderate instances, the new algorithm explores fewer nodes to reach optimality.
In most of these cases, this node reduction translates to smaller execution times, except for the case of \texttt{blacksburg}.
For large instances, node exploration appears important for finding new feasible solutions.
Thus, the dramatic reduction in the number of nodes explored by the new algorithm has a negative impact on the optimality gap reached by the time limit.
Generally, the new algorithm appears useful for instances where optimality can be proven quickly but tends to suffer on large instances because of the increased master problem size --- a topic of future work.

Table \ref{tab:best-bounds} compares the results from the water network design literature with the best results obtained in this study, including where our implementation of \cite{raghunathan2013global} outperformed the new algorithm.
Specifically, the best objective lower and upper bounds are compared.
The table depicts many new bounds discovered using the algorithms herein, especially on outstanding instances.
(For the three \emph{large} cases, our implementation of Raghunathan's algorithm, not the new algorithm, discovered the bounds shown. The exception is for our upper bound on \texttt{pescara}, which the new algorithm discovered.)
Note that here, literature solutions are not differentiated between those obtained heuristically and those obtained via algorithms that can provide lower bounds (e.g., Raghunathan's algorithm).
The slight change in the optimal solution of \texttt{hanoi} could be a typographical error of \cite{raghunathan2013global} or a small difference in our definition of resistance.

\subsection{Comparison of Algorithms on Instances with Scaled Demands}
\label{subsec:scaling}
This section compares the efficacy of the two algorithms on a set of moderate instances extending those described in Section \ref{subsec:experimental_setup}.
For each junction in each network, the original demand was scaled by a factor between $0.5$ and $1.5$ in steps of $0.05$, generating $21$ instances per network and producing $105$ instances.
For \texttt{hanoi}, instances with scalings greater than one were infeasible, reducing the final set to $95$ instances.
For each instance, the time required to prove optimality was measured.

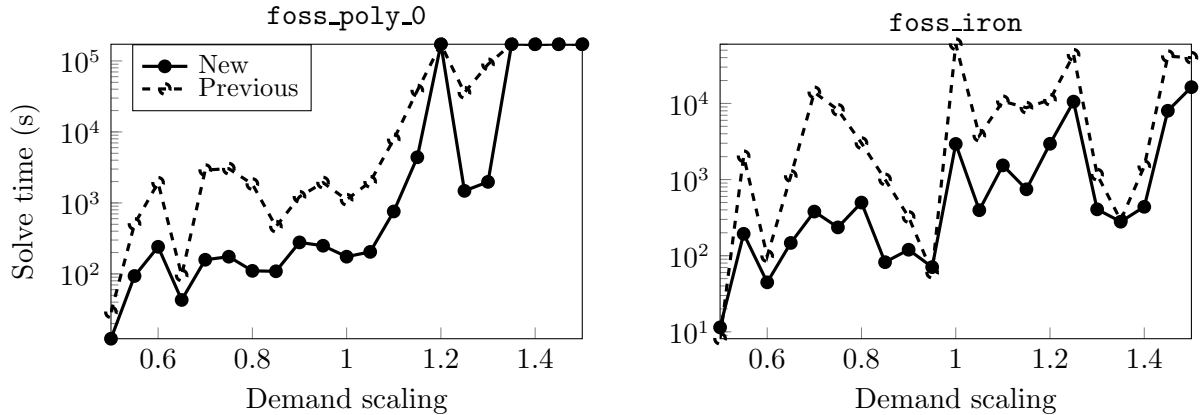
\begin{figure}[!ht]
    \centering
    \begin{subfigure}[b]{0.49\textwidth}
        \centering
        \begin{tikzpicture}[]
	\centering
	\begin{semilogyaxis}[legend cell align=left,enlargelimits=false,xtick pos=left,
                         ytick pos=left,height=5.5cm,ylabel=Solve time (s),xlabel=Demand scaling,
	                     width=0.97\textwidth,legend style={at={(0.045, 1.00)},anchor=north west},legend columns=1,
                         title=\texttt{foss\_poly\_0},
                         title style={yshift=-1.0ex},
				         legend style={font=\small,column sep=1.5pt, row sep=-4.0pt}]
		\pgfplotstableread[col sep = comma]{computational_experiments/data/scalar/NETWORK.foss_poly_0-SIAM.false.csv}{\tasseff};
		\addplot[very thick, mark=*, opacity=1.00] table [x = scalar, y = time_elapsed]{\tasseff};
		\pgfplotstableread[col sep = comma]{computational_experiments/data/scalar/NETWORK.foss_poly_0-SIAM.true.csv}{\siam};
		\addplot[very thick, dashed, mark=o, opacity=1.00] table [x = scalar, y = time_elapsed]{\siam};
		\legend{New,Previous};
	\end{semilogyaxis}
\end{tikzpicture}
    \end{subfigure}
    \hfill
    \begin{subfigure}[b]{0.49\textwidth}
        \centering
        \begin{tikzpicture}[]
	\centering
	\begin{semilogyaxis}[legend cell align=left,enlargelimits=false,xtick pos=left,
                         ytick pos=left,height=5.5cm,xlabel=Demand scaling,
                         title=\texttt{foss\_iron},
                         title style={yshift=-1.0ex},
	                     width=0.97\textwidth]
		\pgfplotstableread[col sep = comma]{computational_experiments/data/scalar/NETWORK.foss_iron-SIAM.false.csv}{\tasseff};
		\addplot[very thick, mark=*, opacity=1.00] table [x = scalar, y = time_elapsed]{\tasseff};
		\pgfplotstableread[col sep = comma]{computational_experiments/data/scalar/NETWORK.foss_iron-SIAM.true.csv}{\siam};
		\addplot[very thick, dashed, mark=o, opacity=1.00] table [x = scalar, y = time_elapsed]{\siam};
	\end{semilogyaxis}
\end{tikzpicture}
    \end{subfigure}
    \caption{Time (log scale) to reach optimality (if attained) across demand scaling instances for \texttt{foss\_poly\_0} and \texttt{foss\_iron}, comparing our new algorithm (New) with an algorithm similar to \cite{raghunathan2013global} (Previous).}
    \label{fig:foss-scalar}
\end{figure}

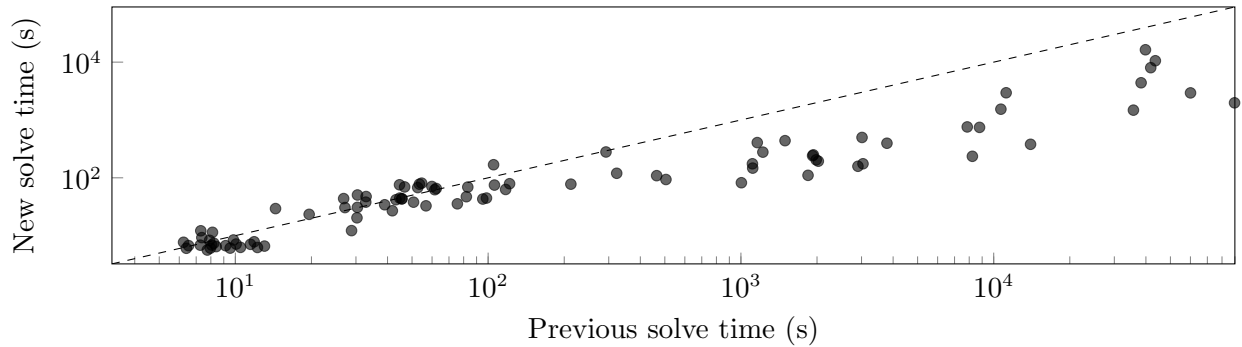
\begin{figure}[t]
     \centering
     \begin{tikzpicture}[]
	\centering
	\begin{loglogaxis}[legend cell align=left,enlargelimits=false,xtick pos=left,
                       ytick pos=left,height=5.0cm,xlabel=Previous solve time (s),
                       ylabel=New solve time (s),width=\textwidth]
		\pgfplotstableread[col sep = comma]{computational_experiments/data/scalar/solve_time_comparison.csv}{\data};
		\addplot[only marks, opacity=0.6] table [x = siam_time, y = tasseff_time]{\data};
        \addplot[domain=3.25:90000.00,samples=100,color=black,dashed]{x};
	\end{loglogaxis}
\end{tikzpicture}
     \vspace{-2em}
     \caption{Time (log scale) to solve demand scaling instances with our new algorithm versus an algorithm similar to \cite{raghunathan2013global} (Previous). Instances where optimality was \emph{not} proven by either are excluded.}
     \label{fig:scalar-time}
\end{figure}

Figure \ref{fig:foss-scalar} compares the times required to reach optimality across demand-scaled versions of \texttt{foss\_poly\_0} and \texttt{foss\_iron}.
These instances display the most dramatic differences between the two algorithms.
For the \texttt{foss\_poly\_0} network, the new algorithm always outperforms the previous algorithm, occasionally by over an order of magnitude.
For demand scalings in the set $\{1.2, 1.35, 1.4, 1.45, 1.5\}$, neither algorithm can solve the design instances within the nearly two day time limit.
The comparison of solve times in the \texttt{foss\_iron} plot implies similar behavior.
For all instances except two ($0.5$ and $0.95$), the new algorithm outperforms the previous algorithm.
Moreover, the new algorithm occasionally outperforms the previous by nearly two orders of magnitude.

Figure \ref{fig:scalar-time} compares the times to reach optimality across all demand-scaled instances.
The figure excludes the five instances where optimality could not be proven.
For instances requiring roughly one hundred seconds or less to reach optimality, the times required by both algorithms are similar, as shown by their placement around the dashed identity line.
For instances requiring roughly one hundred seconds or more, however, the new algorithm always outperforms the previous algorithm, often by one or two orders of magnitude.
This supports Section \ref{subsec:comparison}: for challenging problems where optimality can be proven, our new algorithm outperforms the previous algorithm's implementation.


\section{Concluding Remarks}
\label{section:concluding_remarks}
This study presented the derivation and algorithmic application of a novel, \emph{exact} MICP formulation for the global optimization of potable water distribution network design.
Construction of this problem began with a convex reformulation for network design feasibility and extended the models described by \cite{cherry1951cxvii} and \cite{raghunathan2013global} using nonlinear duality.
Then, using the new MICP formulation as a foundation, the global, linear relaxation-based algorithm of \cite{raghunathan2013global} was augmented to employ novel outer approximation cuts derived from the new formulation.

To measure its efficacy, our global optimization method was compared against the previous state of the art on standard design instances.
Then, new moderately-sized instances were generated by scaling demands throughout the original networks.
The combination of results implies significant speedups in convergence (i.e., one to two orders of magnitude) can be achieved on moderately-sized instances for which optimality can be proven within a modest amount of time (i.e., hours).

This study provides a number of novel and useful contributions to the field of water system design and, more broadly, nonlinear network optimization.
First and foremost, the formulation of a purely \emph{convex} program for determining design feasibility with the inclusion of physical bounds appears to be the first in the literature.
Second and perhaps just as importantly, the exact MICP reformulation of the network design problem establishes a new paradigm for approaching problems of this type.
Third, an algorithm based on the novel features of this MICP is presented and appears capable of efficiently proving global optimality on several challenging instances.
In fact, our computational results close the optimality gap on two outstanding instances (\texttt{foss\_poly\_0} and \texttt{foss\_iron}).

Future work will focus on extending the approaches developed herein.
First, the convex and mixed-integer convex reformulations appear to be immediately applicable to other application areas, including natural gas network expansion planning and crude oil network optimization.
It is also likely that operational problems for nonlinear networks (e.g., pump scheduling optimization in water networks) can be reformulated exactly using methods similar to those developed herein.
Such reformulations may have even greater benefits when considered in these new problem contexts.


\ACKNOWLEDGMENT{%
The authors thank Carleton Coffrin and Kaarthik Sundar of Los Alamos National Laboratory for their expertise in infrastructure optimization and \textsc{JuMP} software development.
This work was conducted under the auspices of the National Nuclear Security Administration of the U.S. Department of Energy at Los Alamos National Laboratory under Contract No. 89233218CNA000001.
Specifically, work at Los Alamos National Laboratory was supported by the Laboratory Directed Research and Development program under the project \emph{Adaptation Science for Complex Natural-engineered Systems} (20180033DR).
}


\section*{Author Biographies}
\paragraph{Byron Tasseff} is a scientist at Los Alamos National Laboratory (LANL) in the Information Systems and Modeling group.
He earned his M.S. in Industrial and Operations Engineering from the University of Michigan in 2018, where he is also a current Ph.D. candidate.
His research involves developing optimization techniques for problems involving fluids and critical infrastructure.

\paragraph{Russell Bent} is a scientist at LANL in the Applied Mathematics and Plasma Physics group.
He is the principal or co-principal investigator for DOE projects in these areas, with focuses on improving robustness of power systems, increasing resilience of infrastructure networks, modeling interdependencies between systems, managing disasters that impact critical infrastructure, modeling smart grid and microgrid technologies, and developing methods for stochastic optimization.

\paragraph{Marina A. Epelman} received her Ph.D. in Operations Research from Massachusetts Institute of Technology. She is currently a professor of Industrial and Operations Engineering at the University of Michigan, with research interests ranging from optimization theory to applications in healthcare and scheduling, among others.

\paragraph{Donatella Pasqualini} is a scientist at LANL in the Information Systems and Modeling group. She  obtained a Ph.D. in Physics from the University of Trento, Italy. Her research focuses on the co-evolution of natural and critical infrastructure systems under climatic changes. 

\paragraph{Pascal Van Hentenryck} is the A. Russell Chandler III Chair and Professor in Industrial and Systems Engineering at the Georgia Institute of Technology, as well as the Associate Chair for Innovation and Entrepreneurship. His current research primarily focuses on large-scale optimization and machine learning applied to energy systems and mobility. 

\paragraph{Story}
This work began as part of a LANL project to develop better climate adaptation strategies for critical infrastructure.
Early work focused on applying techniques from natural gas network optimization to water network design.
After discovering the convex feasibility-testing method of \citet{raghunathan2013global}, work shifted toward developing novel feasibility cuts for this problem.
Pursuing this further ultimately illuminated a subtle property based on duality, which naturally led to the reformulations presented herein.
These reformulations are the first of their kind and have potential applications to other important network problems constrained by nonconvex relationships.


\bibliographystyle{bibliography/informs2014}
\bibliography{bibliography/bibliography.bib}


\begin{APPENDICES}
\section{Derivation of (P(\emph{r}))'s Dual}
\label{section:appendix-dual-derivation}
A straightforward method to derive the dual problem of $(\textnormal{P}(r))$ is via Lagrangian duality, i.e.,
\begin{equation}
    \max_{h} \min_{q^{\pm} \geq 0} \mathcal{L}(q^{+}, q^{-}, h) = \max_{h} g(h) \label{eqn:lagrangian-duality},
\end{equation}
where $\mathcal{L}$ is the Lagrangian of $(\textnormal{P}(r))$ with dual variables $h$ (for flow conservation constraints), and $g(h)$ is the Lagrangian dual function (to later be maximized).
Following $(\textnormal{P}(r))$, its Lagrangian is written as
\begin{equation}
\begin{gathered}
    \mathcal{L}(q^{+}, q^{-}, h) := -\sum_{i \in \mathcal{J}} h_{i} d_{i}
    + \mathlarger{\sum}_{\mathclap{a := (i, j) \in \mathcal{A} : i \in \mathcal{S}}} \left[\frac{L_{a} r_{a}}{1 + \alpha} (q_{a}^{+})^{1 + \alpha} - (h_{i}^{s} - h_{j}) q_{a}^{+} \right]\\
    + \mathlarger{\sum}_{\mathclap{a := (i, j) \in \mathcal{A} : i \in \mathcal{S}}} \left[\frac{L_{a} r_{a}}{1 + \alpha} (q_{a}^{-})^{1 + \alpha} + (h_{i}^{s} - h_{j}) q_{a}^{-} \right]
    + \mathlarger{\sum}_{\mathclap{a := (i, j) \in \mathcal{A} : i \in \mathcal{J}}} \left[\frac{L_{a} r_{a}}{1 + \alpha} (q_{a}^{+})^{1 + \alpha} - (h_{i} - h_{j}) q_{a}^{+} \right]\\
    + \mathlarger{\sum}_{\mathclap{a := (i, j) \in \mathcal{A} : i \in \mathcal{J}}} \left[\frac{L_{a} r_{a}}{1 + \alpha} (q_{a}^{-})^{1 + \alpha} + (h_{i} - h_{j}) q_{a}^{-} \right].
\end{gathered}
\label{eqn:lagrangian}
\end{equation}

Equation \eqref{eqn:lagrangian} is highly separable in $q_{a}^{\pm}$.
As such, minimization over $q^{\pm}$ in Equation \eqref{eqn:lagrangian-duality} is straightforward.
To derive $g(h)$, it suffices to minimize each component of the second through fifth sums over their corresponding $q_{a}^{\pm}$ while imposing nonnegativity on $q_{a}^{\pm}$.
Note that all terms are of the form $\frac{b}{1 + \alpha} \left(q_{a}^{\pm}\right)^{1 + \alpha} + t q_{a}^{\pm}$, where $b > 0$ and the sign of $t$ is unknown.
There are two possibilities: if $t \geq 0$, the component is nondecreasing in $q_{a}^{\pm}$ over $q_{a}^{\pm} \geq 0$, which implies its minimum is attained at $q_{a}^{\pm} = 0$.
Otherwise, if $t < 0$, the function is decreasing at $q_{a}^{\pm} = 0$, attains its minimum, then starts increasing.
This minimum is attained at the stationary point $\hat{q}_{a}^{\pm} = \sqrt[\alpha]{-\frac{t}{b}}$.
In this case, after reduction, the minimum value of the corresponding component is thus
\begin{equation}
    \left(\frac{b}{1 + \alpha} - b\right) \left(-\frac{t}{b}\right)^{1 + \frac{1}{\alpha}} = \frac{-\alpha}{1 + \alpha} \frac{\left(-t\right)^{1 + \frac{1}{\alpha}}}{\sqrt[\alpha]{b}} \label{eqn:dual-min}.
\end{equation}

Next, note that the second and third, as well as the fourth and fifth terms of the sums in Equation \eqref{eqn:lagrangian} can be paired such that the $b$ coefficients of each term are the same, while the $t$ coefficients are opposite.
That is, one term (with nonnegative $t$) has a minimum at zero, while the other has a minimum equal to the right-hand side of Equation \eqref{eqn:dual-min}.
Since the sign of $t$ is unknown, $\lvert t \rvert$ is thus used instead to write $g(h)$ as
\begin{equation}
    g(h) = -\sum_{i \in \mathcal{J}} h_{i} d_{i} - \mathlarger{\sum}_{\mathclap{a := (i, j) \in \mathcal{A} : i \in \mathcal{S}}} \frac{\alpha}{1 + \alpha} \frac{\lvert h_{i}^{s} - h_{j} \rvert^{1 + \frac{1}{\alpha}}}{\sqrt[\alpha]{L_{a} r_{a}}} - \mathlarger{\sum}_{\mathclap{a := (i, j) \in \mathcal{A} : i \in \mathcal{J}}} \frac{\alpha}{1 + \alpha} \frac{\lvert h_{i} - h_{j} \rvert^{1 + \frac{1}{\alpha}}}{\sqrt[\alpha]{L_{a} r_{a}}}.
\end{equation}
Using a standard rewriting of the absolute value terms, the dual problem becomes equivalent to $(\textnormal{D}(r))$.

\section{Physical Interpretation of Strong Duality}
\label{section:appendix-strong-duality}
Four sums appear in the objectives of $(\textnormal{P}(r))$ and $(\textnormal{D}(r))$, each having a unique physical connotation.
This implies a physical meaning of the strong duality constraint.
To begin, let this constraint be expanded as
\begin{equation}
    f_{P}(q) - f_{D}(h) = f_{1}(q) - f_{2}(q) + f_{3}(\Delta h) + f_{4}(h) \leq 0.
\end{equation}
Consider the first summation,
\begin{equation}
    f_{1}(q) = \sum_{a \in \mathcal{A}} \frac{L_{a} r_{a}}{1 + \alpha} \left[(q_{a}^{+})^{1 + \alpha} + (q_{a}^{-})^{1 + \alpha}\right].
\end{equation}
The terms involved are similar to those appearing in the head loss relationships, where $L_{a} r_{a} (q_{a}^{\pm})^{\alpha}$ (conventionally in units of length) can be interpreted as the energy per unit weight of fluid lost to friction between the moving water and the interior wall of the pipe.
Thus, up to a multiplicative constant, $L_{a} r_{a} (q_{a}^{\pm})^{\alpha} q_{a}^{\pm}$ can be interpreted as the rate of heat transference (i.e., power) between the volume of water and the interior wall of the pipe.
This sum can then be interpreted as all power losses from friction of the pipe walls.

Next, consider the second sum appearing in the primal problem $(\textnormal{P}(r))$,
\begin{equation}
    f_{2}(q) = \sum_{i \in \mathcal{S}} h_{i}^{s} \sum_{a \in \delta^{+}_{i}} q_{a} \label{eqn:power-gen}.
\end{equation}
Here, each head $h_{i}^{s}$ is static and can be viewed as the amount of energy per unit weight of water available for extraction from the reservoir.
Thus, its product with the reservoir's (outgoing) flow can be interpreted, again up to a multiplicative constant, as the power generated by the reservoir.
The sum of all contributions in Equation \eqref{eqn:power-gen} is thus proportional to the power supplied to the water network.

Next, consider the first sum appearing in $(\textnormal{D}(r))$, that is,
\begin{equation}
    f_{3}(\Delta h) = \frac{\alpha}{1 + \alpha} \sum_{a \in \mathcal{A}} \frac{1}{\sqrt[\alpha]{L_{a} r_{a}}} \left[(\Delta h^{+}_{a})^{1 + \frac{1}{\alpha}} + (\Delta h^{-}_{a})^{1 + \frac{1}{\alpha}}\right].
\end{equation}
Here, each $\Delta h_{a}^{\pm}$ denotes the difference in energy per unit weight between adjacent nodes.
Each term thus represents, up to a multiplicative constant, the power loss along the pipe in the form of a head differential (i.e., \emph{not} losses to heat from friction).
The sum denotes the total loss in \emph{usable} power across the network.

Finally, consider the second sum appearing in $(\textnormal{D}(r))$, that is,
\begin{equation}
    f_{4}(h) = \sum_{i \in \mathcal{J}} h_{i} d_{i}.
\end{equation}
Here, each demand is fixed, while the energy per unit weight $h_{i}$ at each junction can vary.
Using an argument similar to that of Equation \eqref{eqn:power-gen}, this sum denotes, up to a constant, the power demanded across all junctions.

The strong duality constraint can then be thought of as encoding
\begin{equation}
    \textnormal{(frictional loss)} + \textnormal{(realized loss)} + \textnormal{(demand)} \leq \textnormal{(generation)},
\end{equation}
which implies the conservation of power, with an inequality replacing the traditional equality.

\section{Modified Global Optimization Algorithm}
\label{section:appendix-algorithm}
\begin{algorithm}[!htp]
    \caption{LP/NLP-BB algorithm for the global optimization of $(\textnormal{MINLP})$/$(\textnormal{MICP-E})$.}
    \label{alg:global}
    \begin{algorithmic}[1]
        \State{$r^{*} \gets \textbf{InitialSoln}$ (Algorithm 4 of \citep{raghunathan2013global}); $\eta^{*} \gets \sum_{a \in \mathcal{A}} L_{a} c_{ar_{a}^{*}}$; $\tilde{\mathcal{X}} \gets \emptyset$}. \label{line:mip-init-soln}
        \State{\color{palette2}$(\hat{q}, \hat{h}, \Delta \hat{h}, \hat{x}, \hat{y})$ $\gets$ Solve (MICP-E) with $0 \leq x \leq 1$ using a nonlinear solver.} \label{line:mip-micp-relaxed}
        \State{$r_{a}^{\pm} \gets \argmax_{p \in \mathcal{R}_{a}} \left\{L_{a} p (\hat{q}^{\pm}_{ap})^{\alpha}\right\}, \, \forall a \in \mathcal{A}$; $\tilde{\mathcal{Q}}_{ar_{a}}^{\pm} \gets \{\hat{q}_{ar_{a}}^{\pm}\}, \, \forall a \in \mathcal{A}$.} \label{line:mip-initial-oa-cuts}
        \State{\color{palette2}$r_{a}^{\pm} \gets \argmax_{p \in \mathcal{R}_{a}} \left\{p (\hat{q}^{\pm}_{ap})^{1 + \alpha}\right\}, \, \forall a \in \mathcal{A}$; $\tilde{\mathcal{Q}}_{ar_{a}}^{\textrm{NL}^{\pm}} \gets \{\hat{q}_{ar_{a}}^{\pm}\}, \, \forall a \in \mathcal{A}$.} \label{line:mip-initial-q-nl-cuts}
        \State{\color{palette2}$r_{a}^{\pm} \gets \argmax_{p \in \mathcal{R}_{a}} \left\{\frac{1}{\sqrt[\alpha]{p}} (\Delta \hat{h}_{ap}^{\pm})^{1 + \frac{1}{\alpha}}\right\}, \, \forall a \in \mathcal{A}$; $\mathcal{H}_{ar_{a}}^{\textrm{NL}^{\pm}} \gets \{\Delta \hat{h}_{ar_{a}}^{\pm}\}, \, \forall a \in \mathcal{A}$.} \label{line:mip-initial-dh-nl-cuts}
        \State{\color{palette2}Add $\sum_{a \in \mathcal{A}} L_{a} \sum_{p \in \mathcal{R}_{a}} c_{ap} x_{ap} \geq \sum_{a \in \mathcal{A}} L_{a} \sum_{p \in \mathcal{R}_{a}} c_{ap} \hat{x}_{ap}$ to (MIP-E).} \label{line:mip-lb}
		  \While{$(\textnormal{MIP-ER})$ termination criteria is not satisfied} \label{line:term-criteria}
            \State{$(\hat{q}, \hat{h}, \Delta \hat{h}, \hat{x}, \hat{y})$ $\gets$ Solve the current nodal linear subproblem of $(\textnormal{MIP-ER})$.} \label{line:solve-nodal}
            \State{$\hat{\eta} \gets \sum_{a \in \mathcal{A}} L_{a} \sum_{p \in \mathcal{R}_{a}} c_{ap} \hat{x}_{ap}$.} \label{line:set-objective}
            \If{$\hat{x}_{ap} \in \mathbb{B}, ~ \forall a \in \mathcal{A}, ~ \forall p \in \mathcal{R}_{a}$} \label{line:check-integrality}
                \State{$r_{a} \in \{p \in \mathcal{R}_{a} : \hat{x}_{ap} = 1\}, ~ \forall a \in \mathcal{A}$.} \label{line:get-resistance-sol}
                \State{$(\hat{q}, \hat{h}) \gets$ Solve $(\textnormal{P}(r))$.} \label{line:solve-p-1}
                \If{$\underline{q} \leq \hat{q} \leq \overline{q}$ and $\underline{h} \leq \hat{h} \leq \overline{h}$} \label{line:check-bounds}
                    \State{$r^{*} \gets r$; $\eta^{*} \gets \sum_{a \in \mathcal{A}} L_{a} c_{ar_{a}^{*}}$.} \label{line:mip-set-new-best-sol-1}
                    \State{$\mathcal{\mathcal{Q}}_{ar_{a}}^{\pm} \gets \tilde{\mathcal{Q}}_{ar_{a}}^{\pm} \cup \{\pm\hat{q}_{a}\}, \, \forall a \in \mathcal{A} : \pm \hat{q}_{a} > 0$.} \label{line:mip-oa-1}
                    \State{\color{palette2}$\mathcal{Q}_{ar_{a}}^{\textrm{NL}\pm} \gets \tilde{\mathcal{Q}}_{ar_{a}}^{\textrm{NL}\pm} \cup \{\pm\hat{q}_{a}\}, \, \forall a \in \mathcal{A} : \pm \hat{q}_{a} > 0$.} \label{line:mip-oa-2}
                    \State{\color{palette2}$\tilde{\mathcal{H}}_{ar_{a}}^{\textrm{NL}\pm} \gets \tilde{\mathcal{H}}_{ar_{a}}^{\textrm{NL}\pm} \cup \{\pm\Delta \hat{h}_{a}\}, \, \forall a \in \mathcal{A} : \pm \Delta \hat{h}_{a} > 0$.} \label{line:mip-oa-3}
                \Else
                    \State{\color{palette2}$\tilde{\mathcal{X}} \gets \tilde{\mathcal{X}} \cup \{\hat{x}\}$ (i.e., add a feasibility cut that removes $\hat{x}$).} \label{line:mip-feas-cut}
                    \State{(repaired, $r$) $\gets$ \textbf{Repair}$(r, n, \eta^{*})$ (Algorithm 3 of \citep{raghunathan2013global}).} \label{line:mip-repair-1}
                    \If{repaired} \label{line:mip-if-repair-1}
                        \State{$r^{*} \gets r$; $\eta^{*} \gets \sum_{a \in \mathcal{A}} L_{a} c_{ar_{a}^{*}}$.} \label{line:mip-set-new-best-sol-2}
                        \State{$(\hat{q}, \hat{h}) \gets$ Solve $(\textnormal{P}(r))$.} \label{line:mip-solve-p-2}
            \algstore{myalg}
    \end{algorithmic}
\end{algorithm}
\begin{algorithm}[!htp]
    \caption*{... Continuation of Algorithm \ref{alg:global}.}
    \begin{algorithmic}
        \algrestore{myalg}
                        \State{$\tilde{\mathcal{Q}}_{ar_{a}}^{\pm} \gets \tilde{\mathcal{Q}}_{ar_{a}}^{\pm} \cup \{\pm\hat{q}_{a}\}, \, \forall a \in \mathcal{A} : \pm \hat{q}_{a} > 0$.} \label{line:mip-oa-4}
                        \State{\color{palette2}$\tilde{\mathcal{Q}}_{ar_{a}}^{\textrm{NL}\pm} \gets \tilde{\mathcal{Q}}_{ar_{a}}^{\textrm{NL}\pm} \cup \{\pm\hat{q}_{a}\}, \, \forall a \in \mathcal{A} : \pm \hat{q}_{a} > 0$.} \label{line:mip-oa-5}
                        \State{\color{palette2}$\tilde{\mathcal{H}}_{ar_{a}}^{\textrm{NL}\pm} \gets \tilde{\mathcal{H}}_{ar_{a}}^{\textrm{NL}\pm} \cup \{\pm\Delta \hat{h}_{a}\}, \, \forall a \in \mathcal{A} : \pm \Delta \hat{h}_{a} > 0$.} \label{line:mip-oa-6}
                    \EndIf
                \EndIf
		   	\ElsIf{$(\textnormal{node index}) \mod J = 0$} \label{line:check-node-id}
                \State{\color{palette2}\textbf{NodeCutsNew}$(\tilde{\eta}, \hat{\eta}, m, (\hat{q}, \hat{h}, \Delta \hat{h}, \hat{x}, \hat{y}), \textrm{True})$.} \label{line:call-nodecuts-1}
                \State{$r_{a} \gets \argmin\left\{\hat{x}_{ap} \geq \frac{1}{\lvert \mathcal{R}_{a}\rvert} : p \in \mathcal{R}_{a}\right\}, ~ \forall a \in \mathcal{A}$.} \label{line:get-heuristic-r}
                \State{(repaired, $r$) $\gets$ \textbf{Repair}$(r, \textrm{iter}^{\max}, \eta^{*})$ (Algorithm 3 of \citep{raghunathan2013global}).} \label{line:repair-2}
                \If{repaired} \label{line:if-repair-2}
                    \State{$r^{*} \gets r$; $\eta^{*} \gets \sum_{a \in \mathcal{A}} L_{a} c_{ar_{a}^{*}}$}. \label{line:set-repaired-best-sol-heuristic}
                \EndIf \label{line:endif-heuristic}
            \Else
            \State{\color{palette2}\textbf{NodeCutsNew}$(\tilde{\eta}, \hat{\eta}, m, (\hat{q}, \hat{h}, \Delta \hat{h}, \hat{x}, \hat{y}), \textrm{False})$.} \label{line:call-nodecuts-2}
            \EndIf
            \State{$\tilde{\eta} \gets \hat{\eta}$.} \label{line:store-last-objective}
			\EndWhile
    \end{algorithmic}
\end{algorithm}

Algorithm \ref{alg:global} modifies the algorithm of \cite{raghunathan2013global} to define \emph{this} study's global optimization algorithm.
Portions that substantially modify the original algorithm are denoted by blue font.
The algorithm begins in Line \ref{line:mip-init-soln} by heuristically generating a feasible solution via Algorithm 4 of \cite{raghunathan2013global} and storing the corresponding resistances $r^{*}$ and initial objective $\eta^{*}$.
It also initializes the set of encountered infeasible designs $\tilde{\mathcal{X}}$ to the empty set.
Line \ref{line:mip-micp-relaxed} solves the root relaxation of $(\textnormal{MICP-E})$ via a nonlinear programming algorithm (e.g., \textsc{Ipopt} of \cite{wachter2006implementation}).
Lines \ref{line:mip-initial-oa-cuts} through \ref{line:mip-initial-dh-nl-cuts} use the root relaxation's solution to generate initial linear outer approximations.
In each cut, outer approximation points are chosen to coincide with maximal values of the corresponding nonlinear terms.
Finally, since the root relaxation of $(\textnormal{MICP-E})$ provides a lower bound on the optimal objective, Line \ref{line:mip-lb} imposes this bound explicitly on $(\textnormal{MIP-ER})$.

Line \ref{line:term-criteria} begins the search of a MIP solver, where termination criteria comprises a minimal optimality gap or a time limit.
Line \ref{line:solve-nodal} obtains the relaxation solution for the current node in the search tree.
Line \ref{line:set-objective} computes the objective corresponding to this solution.
Line \ref{line:check-integrality} checks if the current branch and bound (BB) node's resistance choice solution is integer.
If so, Line \ref{line:get-resistance-sol} obtains the corresponding active resistance parameters.
Using these resistances, Line \ref{line:solve-p-1} solves $(\textnormal{P}(r))$ to obtain its primal and dual solutions, $\hat{q}$ and $\hat{h}$, respectively.
This unique solution is then compared against the variable bounds of the original problem in Line \ref{line:check-bounds}.

If bounds are satisfied, $\hat{x}$ is feasible for $(\textnormal{MICP-E})$, and in Line \ref{line:mip-set-new-best-sol-1}, the incumbent solution is updated.
If the solution is physically feasible, Lines \ref{line:mip-oa-1} through \ref{line:mip-oa-3} add outer approximations using the corresponding solution of $(\textnormal{P}(r))$.
Otherwise, if the design is \emph{not} feasible, Line \ref{line:mip-feas-cut} adds a traditional combinatorial no-good feasibility cut.
In Line \ref{line:mip-repair-1}, Algorithm 3 of \cite{raghunathan2013global} is called in the attempt to heuristically recover a new feasible incumbent solution by ``repairing'' infeasibilities of the design $\hat{x}$, where $n$ denotes the (fixed) maximum number of repair iterations to be used by the procedure.
If the solution \emph{was} repaired, as indicated in Line \ref{line:mip-if-repair-1}, then a new incumbent was found, which is updated in Line \ref{line:mip-set-new-best-sol-2}.
Using the resistances from this design, Line \ref{line:mip-solve-p-2} solves $(\textnormal{P}(r))$ to obtain the exact physical solution of the network.
Then, in Lines \ref{line:mip-oa-4} through \ref{line:mip-oa-6}, outer approximations are added to $(\textnormal{MIP-ER})$ based on the physical solution corresponding to the design.

If the design solution is not integral, Line \ref{line:check-node-id} serves as a possible entry point for adding outer approximations to $(\textnormal{MIP-ER})$ and heuristically discovering new solutions.
In this case, if the integer index of the BB node is divisible by some integer $J$, these routines are called.
In Line \ref{line:call-nodecuts-1}, Algorithm \ref{alg:nodecuts} is called, which adds outer approximations based on the relaxation solution at the current BB node.
This algorithm is described in Appendix \ref{appendix:nodecuts}.
Then, in Line \ref{line:get-heuristic-r}, a heuristic resistance solution is prepared, which selects from active resistances in the current relaxation solution.
The repair algorithm is then invoked on this (inexpensive but presumably infeasible) network design in Lines \ref{line:repair-2} through \ref{line:set-repaired-best-sol-heuristic}.
For greater detail, this heuristic procedure is developed and elaborated upon by \cite{raghunathan2013global}.
Otherwise, if both the design solution is fractional and the above heuristic is not activated, Algorithm \ref{alg:nodecuts} is called on Line \ref{line:call-nodecuts-2}, which conditionally refines $(\textnormal{MIP-ER})$'s outer approximations using the solution of the continuous relaxation at the current BB node.

Finally, on Line \ref{line:store-last-objective}, the objective at the current BB node, $\hat{\eta}$, is stored as $\tilde{\eta}$ for use within methods of the next node.
Specifically, both the current objective $\hat{\eta}$ and previous objective $\tilde{\eta}$ are used in a conditional step of Algorithm \ref{alg:nodecuts} to determine whether outer approximations should be added to $(\textnormal{MIP-ER})$.

Algorithm \ref{alg:global} is highly similar to the algorithm of \cite{raghunathan2013global} but differs in a few important respects.
First, in Lines \ref{line:mip-initial-oa-cuts} through \ref{line:mip-initial-dh-nl-cuts}, initial outer approximations are based on the root relaxation of $(\textnormal{MICP-E})$, whereas Raghunathan's algorithm uses the root relaxation of $(\textnormal{MICP-R})$.
Second, Raghunathan only applies outer approximations similar in form to Constraints (22), which correspond to convex head loss relaxations.
Algorithm \ref{alg:global} extends this by adding outer approximations of terms appearing in the strong duality Constraint (27).
Lastly, Line \ref{line:call-nodecuts-1} ensures outer approximations will occasionally be added to $(\textnormal{MIP-ER})$, whereas the conditions in Raghunathan's algorithm (and in their \textbf{NodeCuts} algorithm) are more restrictive.

In our study, implementations of the new and previous algorithms required some modifications, as necessitated by feature limitations of \textsc{JuMP}.
First, the depth of a node in the branch and bound tree ($m$ in Algorithm \ref{alg:global}) is used as a parameter in both algorithms but is not easily accessible via \textsc{JuMP}.
Thus, the number of unity-valued components of $\hat{x}$ at a node in the BB tree is used instead.
Additionally, a number of modifications are made to Raghunathan's algorithm to ensure fairer comparison with Algorithm \ref{alg:global}.
The most substantial of these are as follows: (i) initial outer approximation points are taken from the root relaxation of $(\textnormal{MICP-E})$; and (ii) similar to Lines \ref{line:check-node-id} through \ref{line:endif-heuristic} in Algorithm \ref{alg:global}, the heuristic \emph{and} \textbf{NodeCuts} methods are called every $J$ BB nodes.
Aside from these changes, Raghunathan's algorithm is nearly reproduced.

\section{Modified NodeCuts Algorithm}
\label{appendix:nodecuts}
\begin{algorithm}[!ht]
    \caption{$\textbf{NodeCutsNew}(\tilde{\eta}, \hat{\eta}, m, (\hat{q}, \hat{h}, \Delta \hat{h}, \hat{x}, \hat{y}), \color{palette2}\textrm{Force})$}
    \label{alg:nodecuts}
    \begin{algorithmic}[1]
        \If{$\left[\left(\textrm{rand}\left([0,1]\right) \leq \beta_{oa} 2^{-m}\right) \land \left(\left\lvert\frac{\hat{\eta} - \tilde{\eta}}{\tilde{\eta}}\right\rvert \geq K_{oa}\right)\right] \color{palette2}{\lor \; \textrm{Force} = \textrm{True}}$} \label{line:nodecuts-comparison}
			\For{$a \in \mathcal{A}$} \label{line:nodecuts-loop-start}
            \If{$\hat{y}_{a} \geq 0.5$} \label{line:nodecuts-if-positive-flow}
                    \State{$r_{a} \gets \argmax_{p \in \mathcal{R}_{a}} \left\{L_{a} p (\hat{q}^{+}_{ap})^{\alpha}\right\}$} \label{line:nodecuts-get-p-r-1}
                    \State{\textbf{if} $\left\lvert L_{a} r_{a} \left(\hat{q}^{+}_{ar_{a}}\right)^{\alpha} - \Delta \hat{h}_{ar_{a}}^{+}\right\rvert > \epsilon$ \textbf{then} $\tilde{\mathcal{Q}}_{ar_{a}}^{+} \gets \tilde{\mathcal{Q}}_{ar_{a}}^{+} \cup \{\hat{q}_{ar_{a}}^{+}\}$ \textbf{end if}} \label{line:nodecuts-add-p-1}
                    \State{\color{palette2}$r_{a} \gets \argmax_{p \in \mathcal{R}_{a}} \left\{p (\hat{q}^{+}_{ap})^{1 + \alpha}\right\}$} \label{line:nodecuts-get-p-r-2}
                    \State{\color{palette2}\textbf{if} $\left\lvert \frac{r_{a} \left(\hat{q}^{+}_{ar_{a}}\right)^{1 + \alpha}}{1 + \alpha} - \hat{q}_{a}^{\textrm{NL}} \right\rvert > \epsilon$ \textbf{then} $\tilde{\mathcal{Q}}_{ar_{a}}^{\textrm{NL}+} \gets \tilde{\mathcal{Q}}_{ar_{a}}^{\textrm{NL}+} \cup \{\hat{q}^{+}_{ar_{a}}\}$ \textbf{end if}} \label{line:nodecuts-add-p-2}
                    \State{\color{palette2}$r_{a} \gets \argmax_{p \in \mathcal{R}_{a}} \left\{\frac{1}{\sqrt[\alpha]{p}} (\Delta \hat{h}^{+}_{ap})^{1 + \frac{1}{\alpha}}\right\}$} \label{line:nodecuts-get-p-r-3}
                    \State{\color{palette2}\textbf{if} $\left\lvert \frac{(\Delta \hat{h}^{+}_{ar_{a}})^{1 + \frac{1}{\alpha}}}{(1 + \frac{1}{\alpha})\sqrt[\alpha]{r_{a}}} - \Delta \hat{h}_{a}^{\textrm{NL}} \right\rvert > \epsilon$ \textbf{then} $\tilde{\mathcal{H}}_{ar_{a}}^{\textrm{NL}+} \gets \tilde{\mathcal{H}}_{ar_{a}}^{\textrm{NL}+} \cup \{\Delta \hat{h}^{+}_{ar_{a}}\}$ \textbf{end if}} \label{line:nodecuts-add-p-3}
				\Else \label{line:nodecuts-else}
                    \State{$r_{a} \gets \argmax_{p \in \mathcal{R}_{a}} \left\{L_{a} p (\hat{q}^{-}_{ap})^{\alpha}\right\}$} \label{line:nodecuts-get-n-r-1}
                    \State{\textbf{if} $\left\lvert L_{a} r_{a} \left(\hat{q}^{-}_{ar_{a}}\right)^{\alpha} - \Delta \hat{h}_{ar_{a}}^{-}\right\rvert > \epsilon$ \textbf{then} $\tilde{\mathcal{Q}}_{ar_{a}}^{-} \gets \tilde{\mathcal{Q}}_{ar_{a}}^{-} \cup \{\hat{q}_{ar_{a}}^{-}\}$ \textbf{end if}} \label{line:nodecuts-add-n-1}
                    \State{\color{palette2}$r_{a} \gets \argmax_{p \in \mathcal{R}_{a}} \left\{p (\hat{q}^{-}_{ap})^{1 + \alpha}\right\}$} \label{line:nodecuts-get-n-r-2}
                    \State{\color{palette2}\textbf{if} $\left\lvert \frac{r_{a} \left(\hat{q}^{-}_{ar_{a}}\right)^{1 + \alpha}}{1 + \alpha} - \hat{q}_{a}^{\textrm{NL}} \right\rvert > \epsilon$ \textbf{then} $\tilde{\mathcal{Q}}_{ar_{a}}^{\textrm{NL}-} \gets \tilde{\mathcal{Q}}_{ar_{a}}^{\textrm{NL}-} \cup \{\hat{q}^{-}_{ar_{a}}\}$ \textbf{end if}} \label{line:nodecuts-add-n-2}
                    \State{\color{palette2}$r_{a} \gets \argmax_{p \in \mathcal{R}_{a}} \left\{\frac{1}{\sqrt[\alpha]{p}} (\Delta \hat{h}^{-}_{ap})^{1 + \frac{1}{\alpha}}\right\}$} \label{line:nodecuts-get-n-r-3}
                    \State{\color{palette2}\textbf{if} $\left\lvert \frac{(\Delta \hat{h}^{-}_{ar_{a}})^{1 + \frac{1}{\alpha}}}{(1 + \frac{1}{\alpha})\sqrt[\alpha]{r_{a}}} - \Delta \hat{h}_{a}^{\textrm{NL}} \right\rvert > \epsilon$ \textbf{then} $\tilde{\mathcal{H}}_{ar_{a}}^{\textrm{NL}-} \gets \tilde{\mathcal{H}}_{ar_{a}}^{\textrm{NL}-} \cup \{\Delta \hat{h}^{-}_{ar_{a}}\}$ \textbf{end if}} \label{line:nodecuts-add-n-3}
				\EndIf
			\EndFor
		 \EndIf
    \end{algorithmic}
\end{algorithm}

Algorithm \ref{alg:nodecuts} modifies the \textbf{NodeCuts} algorithm of \cite{raghunathan2013global} to define this study's implementation.
Portions that substantially modify the original are denoted by blue font.
The method requires five parameters: $\tilde{\eta}$, the objective value of the previous BB node's relaxation solution; $\hat{\eta}$, the objective value at the current node; $m$, the depth of the current node; $(\hat{q}, \hat{h}, \Delta \hat{h}, \hat{x}, \hat{y})$, the relaxation solution at the current node; and $\textrm{Force}$, a Boolean variable describing whether or not to force the addition of cuts.
Line \ref{line:nodecuts-comparison} begins by checking three conditions.
First, a random number uniformly generated between zero and one is compared against the term $\beta_{oa} 2^{-m}$, where $\beta_{oa}$ is a positive parameter.
This encourages cuts to be added near the top of the search tree.
The second condition, $\left\lvert \frac{\hat{\eta} - \tilde{\eta}}{\tilde{\eta}} \right\rvert \geq K_{oa}$, computes the relative change in the objective between successive nodes and compares it with the positive constant $K_{oa}$.
This encourages cuts to be added when the solution has significantly changed.
Finally, $\textrm{Force}$ can override the prior conditions, as used in Line \ref{line:call-nodecuts-1} of Algorithm \ref{alg:global}.

If conditions have been satisfied, the algorithm proceeds in Line \ref{line:nodecuts-loop-start} by looping over all arcs.
In Line \ref{line:nodecuts-if-positive-flow}, if the relaxed flow solution along an arc is overall positive, cuts for positive flow are potentially added in Lines \ref{line:nodecuts-add-p-1}, \ref{line:nodecuts-add-p-2}, and \ref{line:nodecuts-add-p-3}.
The reference resistances for each of these cuts, $r_{a}$, are independently computed in Lines \ref{line:nodecuts-get-p-r-1}, \ref{line:nodecuts-get-p-r-2}, and \ref{line:nodecuts-get-p-r-3}, where each corresponds to the maximum nonlinear term given the current relaxation solution.
The conditionals on Lines \ref{line:nodecuts-add-p-1}, \ref{line:nodecuts-add-p-2}, and \ref{line:nodecuts-add-p-3} imply these cuts will only be added if the relaxation solution deviates substantially (greater than a violation of some small constant, $\epsilon$) compared to the solution of a corresponding nonlinear formulation.
This limits the number of cuts being added by ensuring that only constraints with significant violations will be linearized.
In Lines \ref{line:nodecuts-get-n-r-1} through \ref{line:nodecuts-add-n-3}, the above process is analogously completed for arcs where flow along that arc is overall negative.
In our study, parameters of Algorithms \ref{alg:global} and \ref{alg:nodecuts} coincided with those used by Raghunathan, namely, $\beta_{oa} = 5$, $J = 500$, $K_{oa} = 10^{-3}$, $n = 50$, and $\epsilon = 10^{-6}$.
\end{APPENDICES}

\end{document}